\documentclass[10pt,a4paper]{article}
\usepackage[a4paper,margin=0.45in]{geometry}

\usepackage{graphicx}
\usepackage{subcaption}
\usepackage{colortbl}
\usepackage{float}
\usepackage{url,lineno,microtype}
\usepackage[english]{babel}
\usepackage{amssymb,amsmath,amsthm,amsfonts}
\usepackage{color,soul}
\usepackage{adjustbox,lipsum} 
\usepackage{esdiff}
\usepackage{physics}
\usepackage[version=3]{mhchem}

\usepackage{siunitx}
\usepackage{wrapfig}
\usepackage{chemformula} 
\usepackage{mathptmx,cuted}
\usepackage{physics}
\usepackage{sectsty}
\usepackage{fancyhdr}
\usepackage{fnpos}
\usepackage{array}
\usepackage{droidsans}
\usepackage{charter}
\usepackage{epstopdf}
\usepackage{esdiff}
\usepackage{adjustbox,lipsum} 
\usepackage{xcolor}
\usepackage[skip=2pt,font=footnotesize]{caption}
\usepackage[super,sort&compress]{natbib} 

\let\oldsqrt\sqrt
\def\sqrt{\mathpalette\DHLhksqrt}
\def\DHLhksqrt#1#2{%
	\setbox0=\hbox{$#1\oldsqrt{#2\,}$}\dimen0=\ht0
	\advance\dimen0-0.2\ht0
	\setbox2=\hbox{\vrule height\ht0 depth -\dimen0}%
	{\box0\lower0.4pt\box2}}
\makeatletter
\usepackage{array}
\newcolumntype{P}[1]{>{\centering\arraybackslash}p{#1}}

\usepackage{xr-hyper}
\usepackage{hyperref}

\usepackage{authblk} 

\title{A Bilayer Cathode Design Procedure for Li ion Batteries Using the Multilayer Doyle-Fuller-Newman Model (M-DFN)}  

\date{}
\author[1,2,5]{E.~C.~Tredenick\thanks{Email: \texttt{Eloise.Tredenick@canberra.edu.au}}}
\author[3]{A.~M.~Boyce}
\author[4]{R.~Drummond}
\author[1,5]{S.~R.~Duncan}
 \affil[1]{\small Department of Engineering Science, University of Oxford, Oxford, OX1 3PJ, UK}
\affil[2]{Faculty of Science and Technology, University of Canberra, Bruce, 2617, Australia}
\affil[3]{School of Mechanical and Materials Engineering, University College Dublin, Dublin, Ireland}
\affil[4]{School of Electrical and Electronic Engineering, University of Sheffield, Sheffield, S1 3JD, UK}
\affil[5]{The Faraday Institution, Quad One,  Harwell Campus, Didcot, OX11 0RA, UK}

\begin{document}	
\onecolumn
\maketitle

\begin{abstract}

Heterogeneities in lithium ion batteries can be significant factors in electrode under utilisation and degradation while charging. Bilayer electrodes have been proposed as a convenient and scalable way to homogenise the electrode response. In this paper, the design of a bilayer cathode for Li-ion batteries composed of separate layers of lithium nickel manganese cobalt oxide Li[Ni$_{0.6}$Mn$_{0.2}$Co$_{0.2}$]O$_{2}$ (NMC622) and lithium iron phosphate \ce{LiFePO4} (LFP) is optimised using the multilayer Doyle-Fuller-Newman (M-DFN) model. Changes to the carbon binder domain, electrolyte volume fraction, and tortuosity provided the greatest control for improving Li-ion charge mobility. The optimised bilayer design was able to charge at 3C between 0-90\% SOC in 18.6 minutes, achieving 4.4~mAh/cm$^2$. Comparing the optimal bilayer to the LFP-only electrode, the bilayer achieved 41\% higher capacity. Through mechanistic physics-based modelling, it was shown that the 3C charging improvement of the optimised bilayer was achieved by enabling a more homogeneous current density distribution through the thickness of the electrode and electrolyte depletion prevention. The findings were confirmed on a high-fidelity X-ray computed tomography (CT) based microstructural model. The results illustrate how modelling can be used to rapidly search novel electrode designs and accelerate the deployment of fast-charging thick electrodes by adapting existing manufacturing processes.

	\end{abstract}
	
{\footnotesize \textit{Lithium-ion, battery, fast charging, optimisation, bilayer, Multilayer Doyle-Fuller-Newman Model, M-DFN}}

\qquad

Transitioning to electric vehicles (EVs) and replacing internal combustion engines to achieve net zero carbon emissions in the 2035 to 2050 time frame\cite{COP26} is an urgent challenge for addressing climate change. To support this transition, there has been growing focus on improving electrical energy storage technologies, with lithium-ion batteries (LIBs) currently being preferred due to their high energy density and efficiency\cite{hadjipaschalis2009overview,diouf2015potential}. Satisfying the needs of the growing market is driving research to further enhance LIBs. In addition to developing new battery chemistries, it is important to consider existing, well-characterised materials that are scalable with secure supply chains, and explore innovative designs and manufacturing techniques\cite{IOPDrummond2022,Tredenick2024multilayer,Chowdhury2021}. Specifically, the use of tailored arrangements of electrode materials to reduce inhomogeneous responses, and deliver more efficient charge storage, has been proposed. In particular, using physics-based modelling\cite{Tredenick2024multilayer}, electrode bilayers with different chemistries in each sub-layer have been shown in model simulations and experiments to homogenise active material utilisation and electrolyte concentration, improve capacity retention\cite{Tredenick2024multilayer}, and slow degradation experimentally\cite{Enpower2020}. These earlier studies considered only a limited number of electrode designs (for example, 50/50 bilayer thickness) but even then, the designs demonstrated the opportunities of mathematical modelling to rapidly screen potential designs and investigate the underlying mechanisms driving the responses\cite{Tredenick2024multilayer}.

\begin{figure} [h!]
	\centering
	\includegraphics[width=0.45\textheight,keepaspectratio]{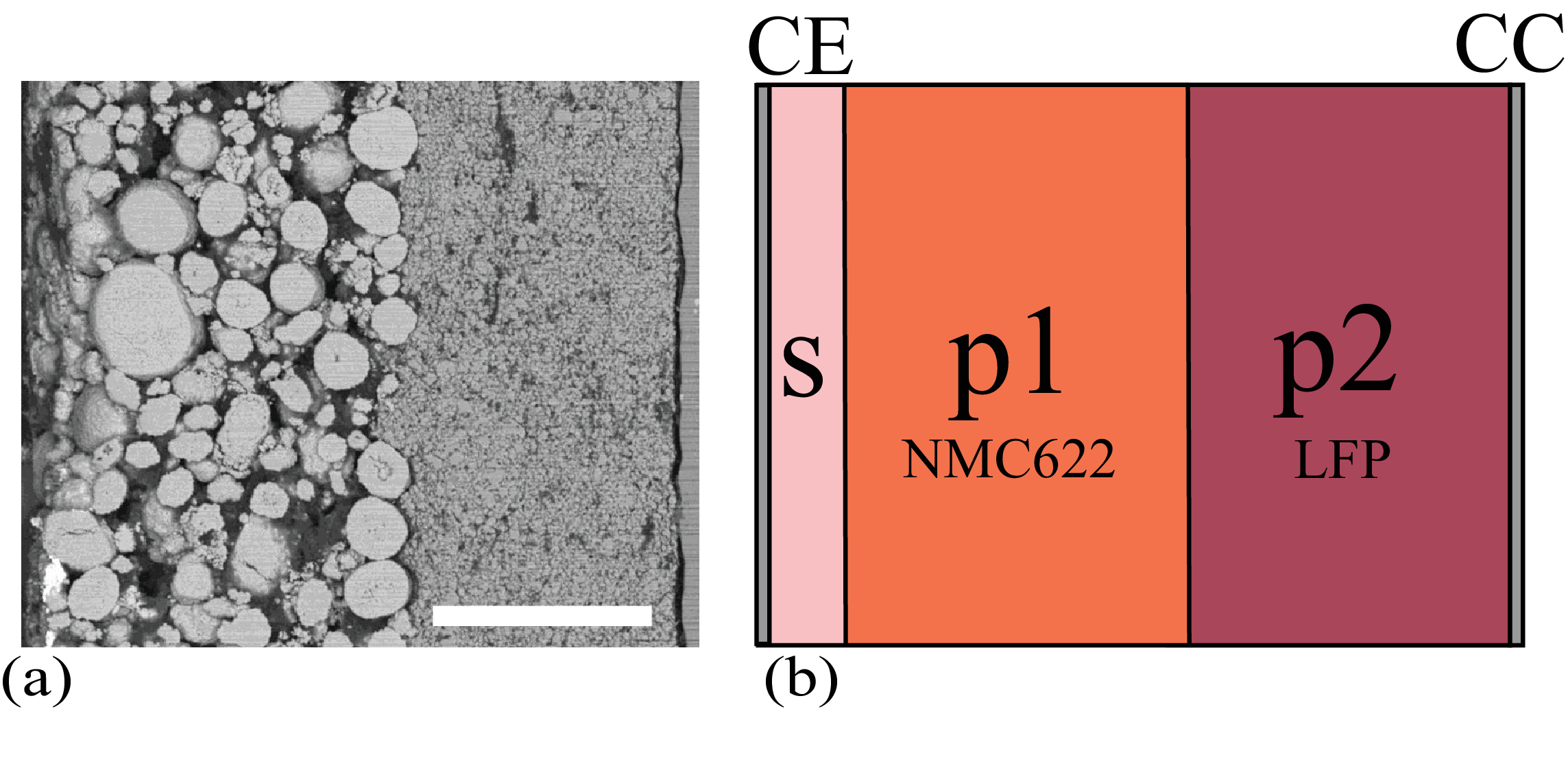}
	\caption{ (a) Plasma focused ion beam/scanning electron microscope (PFIB/SEM) back scattered electron (BSE) image of a bilayer positive electrode for Li ion half-cell comprising a \ce{Li(NiMnCo)O2} (NMC) sub-layer adjacent to the separator, and a \ce{LiFePO4} (LFP) sub-layer next to the Al current collector (CC), with a scale bar of 50~$\mu$m. (b) The M-DFN model domain, comprising Li counter electrode (CE), separator (s), positive electrode sub-layer p1 (NMC in this case), positive electrode sub-layer p2 (LFP in this case), and current collector (CC). There are also spherical particles at each point through the thickness of the electrode for the additional pseudo-dimension (not shown). The model includes three distinct compartments including (i) separator, (ii) positive one, and (iii) positive two. Reproduced with permission\cite{Tredenick2024multilayer}.}
	\label{fig:DomainhalfSEM}  
\end{figure}

Mathematical modelling is becoming increasingly recognised as a cost-effective and fast method to explore and optimise complex battery electrode arrangements\cite{Tredenick2024multilayer,Tredenick2025CTvsDFN,Chowdhury2021Sim,knehr2023material, ge2023numerical,quarti2023trade,de2013model,du2017understanding}. Using the `Multilayer Doyle-Fuller-Newman' (M-DFN) model\cite{Tredenick2024multilayer}, we found that compared with an electrode consisting of a single layer, a bilayer consisting of separate layers of lithium nickel manganese cobalt oxide Li[Ni$_{0.6}$Co$_{0.2}$Mn$_{0.2}$]O$_{2}$ (NMC622) and lithium iron phosphate \ce{LiFePO4} (LFP). Higher capacity is achieved by placing high capacity NMC622 close to the separator, so that the lower capacity LFP next to the current collector receives lower currents but higher active particle utilisation. This is due to their open circuit potential profiles operating in different voltage windows\cite{Tredenick2024multilayer}. Fig.~\ref{fig:DomainhalfSEM} shows the bilayer arrangement including the modelling domain. This existing\cite{Tredenick2024multilayer} bilayer design was not formally optimised, and it was predicted that further performance improvements could be achieved. To address this gap, the focus of this paper is on developing a model-led framework to optimise the design of bilayer cathodes\cite{Tredenick2024multilayer}. We include an X-ray computed tomography (CT) based microstructural model focusing on the default and optimal cases to further validate the M-DFN model results, which is expanding upon the previously developed by the authors\cite{Tredenick2025CTvsDFN}.

The model-led analysis on bilayer considered in this work builds upon existing studies. For example, a multi-objective function and a Doyle-Fuller-Newman (DFN) model has been used to maximise the normalised energy density and minimise the maximum temperature in a battery pack\cite{astaneh2022multiphysics}. The optimisation of four design parameters (including thickness and porosity of the anode and cathode) to maximise the energy density of Li ion batteries using a reformulated DFN model has also been investigated\cite{de2013model}. The balance between a battery's energy density and fast-charging performance was investigated using a DFN model including a range of charge protocols and thicknesses to find an optimal EV battery\cite{quarti2023trade}. The optimisation of microstructural parameters to improve energy density has been investigated using generative artificial intelligence including a Bayesian optimisation loop \cite{kench2024li}. A DFN-type model including double layer capacitance, film resistance and side reactions has been used to determine the optimal energy density of a LFP and graphite cell\cite{srinivasan2004design}. Comparing our work to previous studies, we investigate the optimisation of an experimentally validated model\cite{Tredenick2024multilayer} of a bilayer cathode composed of LFP and NMC sub-layers through a process that considers the electrode design parameters that can be controlled during manufacturing. By examining the optimal bilayer compared to the original bilayer and conventional electrodes, we study the underlying physical processes such as electrolyte concentration and potential, and particle surface reactions to explain the higher capacity. We also expand the previously validated CT model\cite{Tredenick2025CTvsDFN} based on conventional LFP and NMC electrodes to include the bilayer.

\section*{Methodology }

We investigate the potential of the `Multilayer Doyle-Fuller-Newman' (M-DFN) pseudo-2-dimensional (P2D) model for optimising the design of bilayer lithium-ion battery cathodes\cite{Tredenick2024multilayer,Wheeler2023}. The M-DFN model was validated against experimental data across a range of configurations and C rates\cite{Tredenick2024multilayer,Wheeler2023} and is used here to increase the capacity of the cathode at high C rates. With a view to supporting EV fast charging, the proposed M-DFN model derived designs are optimised to maximise 3C (20 minutes) charging capacity. We aim to design the thickest high energy density electrode that achieves increased active material utilisation, without sacrificing capacity retention.

We first conduct a sensitivity analysis to consider a wide range of model parameters of the M-DFN model for 3C charging, with an equivalent specific capacity at 0.05C carefully maintained across cases. We then focus on parameters that increase capacity. Secondly, we combine several parameters and generate a new candidate optimal design. Thirdly, we take the candidate optimal design and investigate the effect of electrode thickness and sub-layer thickness ratio on capacity. The cost function being optimised is the normalised capacity retention (achieved capacity normalised by specific capacity at 0.05C), where the specific capacities at 0.05C are equivalent. When increasing electrode thickness, specific capacities cannot be equivalent, so we consider both capacity and normalised capacity. As we target relatively high energy density cells for EV applications, only small ($<$10\%) changes to porosity are considered.

The design procedure is applied by calculating the new set of parameters and evaluating how they impact the model simulation results. We use the terminology \textit{`benchmark comparison'} when evaluating the performance designs with equivalent specific areal capacity and C rate. As the benchmark comparison requires the specific capacity at 0.05C to be equivalent between cases, the parameters for each case were calculated prior to simulation to enable a fair comparison. The bilayer electrode arrangement that is the focus of this work is termed `NMC:LFP:CC', with a LFP sub-layer adjacent to the metallic current collector (CC) and a NMC622 sub-layer adjacent to the separator and lithium metal counter electrode (CE). `NMC-only' and `LFP-only' are conventional non-layered single active material electrodes. Additional configurations including swapped LFP:NMC:CC and blended/mixed NMC and LFP (NMC+LFP) are not  studied in detail here as they were deemed inferior in previous modelling and experimental studies\cite{Wheeler2023}. The novelty of this work is using the normalised capacity as the cost function, identifying key optimisation parameters to increase capacity and identifying a new optimal bilayer candidate that can be fabricated in the future. We consider a wide range of parameters and carefully maintain the specific capacity at 0.05C to investigate the impact of parameter variations at a 3C charge, along with testing porosity and tortuosity at a range of C rates. We also include an X-ray computed tomography (CT) based microstructural model focusing on the default and optimal cases to further validate the M-DFN model results. All the modelling results including the optimisation were found primarily using the M-DFN model\cite{Tredenick2024multilayer}, and then the CT model was utilised after the M-DFN optimisation routine was concluded to reproduce and validate the M-DFN model optimal case.

\subsection*{Multilayer Doyle-Fuller-Newman Model (M-DFN) }

We solve the pseudo-two dimensional M-DFN, as described in our earlier work\cite{Tredenick2024multilayer}. Briefly, the M-DFN model (Fig.~\ref{fig:DomainhalfSEM}), has two distinct layers of different chemistries including different particles, OCP and diffusivities, the electrolyte concentration and potential are functions of the Li ions in the electrolyte, and includes carbon binder. The model parameters for the bilayer are shown in Table~\ref{VariablesModel} and \ref{ParametersModel}, along with new thicknesses and specific capacities in Table~\ref{BilNMCLFP}. The parameters for the conventional cells are identical to the previous model\cite{Tredenick2024multilayer} unless stated later in Tables \ref{BilNMCLFP} and \ref{3CResults}. The MATLAB code is the same as the previous model\cite{Tredenick2024multilayer}, except for the new parameters, and available online as open-source code for the bilayer \url{https://github.com/EloiseTredenick/M-DFN_Matlab_Code_Li_Ion_Batteries} and uniform cells \url{https://github.com/EloiseTredenick/DFN-P2D-Uniform-Matlab-NMC-LFP-LiIon-Batteries}.

 \subsection*{CT Image-Based Model }

The CT image-based model is constructed using NMC and LFP images from a previous study\cite{Tredenick2025CTvsDFN} and solved using the finite element method. Simpleware ScaniP (Mountain View, CA, USA) is used to combine and crop the NMC and LFP images to the appropriate dimensions. Simpleware is then used to discretise the segmented images giving linear tetrahedral meshes with the 2.4 million mesh elements and 1.5 million degrees of freedom for the 50:50 case, and 3 million mesh elements and 1.7 degrees of freedom for the optimal case. A mesh sensitivity test was carried out and the simulation conditions were found to be stable. The electrochemical and transport as outlined in our previous study\cite{Tredenick2025CTvsDFN} formed the basis for the current study and are not shown here for brevity. The equations are solved in COMSOL\textsuperscript \textregistered  MultiPhysics (v6.1, Sweden). Time stepping is handled using 2nd order backward Euler differentiation. The active material is resolved in both NMC and LFP images and the domain containing pores and carbon binder domain are homogenised\cite{Tredenick2025CTvsDFN}. This region is modelled using a homogenised approach, similar to the P2D DFN model, where an averaged value of the volume fraction of CBD and electrolyte is used throughout the domain. The volume fraction of CBD, along with all other material parameters are outlined in Table \ref{optimaltable4}. The CT model was utilised only in the case of the default bilayer and optimal bilayer, while all other results and optimisation were primarily conducted using the P2D M-DFN model, for computational efficiency. 

\begin{table} [h!]
	\scriptsize
	\centering
	\caption[Model Parameters]{Model parameters for the NMC:LFP:CC M-DFN bilayer electrode (44~$\mu$m:44~$\mu$m, total electrode thickness 88~$\mu$m) for Fig.~\ref{fig:newcapV3C}, where p1 and p2 are NMC and LFP, respectively and are identical to the previous work\cite{Tredenick2024multilayer}.}	 	\label{ParametersModel}	
	\begin{tabular}{    p{1.3cm}  p{7.5cm}  p{2cm}  p{2.5cm}  p{2.2cm} } 
		\hline
		Parameter & Description  & Unit  & Value & Source   \\ \hline 
		$A$& Electrode cross-sectional area & m$^2$ & $1.54\times 10^{-4}$ & Experiment   \\ 	 
		$b_{\scriptscriptstyle \text{p1}}$& Bruggeman tortuosity factor & - & 1.6 &  ${}$\cite{Korotkin2021} \\  
		$b_{\scriptscriptstyle \text{p2}}$& & - & 2.1&  ${}$\cite{Tredenick2024multilayer}   \\ 
		$b_{\scriptscriptstyle \text{s}}$& & - & 1.5 & ${}$\cite{kirk2021physical2}   \\ 
		$c_{\scriptscriptstyle \text{e0}}$&Initial concentration of Li ions in electrolyte & mol/m$^3$ & 1000 &Experiment  \\ 
		$ c_{\scriptscriptstyle \text{s,p1,0}}$& Initial concentration of Li ions in solid particles, charge and discharge &  mol/m$^3$ & 44868, 13366  &  OCP function\cite{Tredenick2024multilayer} \\ 
		$ c_{\scriptscriptstyle \text{s,p2,0}}$&  &  mol/m$^3$ & 22751,   29 &  OCP function\cite{Tredenick2024multilayer}  \\ 
		$ c^{\scriptscriptstyle \text{max}}_{\scriptscriptstyle \text{s,p1}}$& Maximum concentration of Li ions in solid particles &  mol/m$^3$ & 48700  & ${}$\cite{Xu2019}    \\ 
		$ c^{\scriptscriptstyle \text{max}}_{\scriptscriptstyle \text{s,p2}}$&   &  mol/m$^3$ & 22806  & ${}$\cite{Jokar2018mesoscopic}   \\ 	
		$D_{\scriptscriptstyle \text{s,p1}}$&Diffusivity of Li ions in solid & m$^2$/s & $4\times 10^{-14}$ & ${}$\cite{Noh2013comparison}   \\ 
		$D_{\scriptscriptstyle \text{s,p2}}$&  & m$^2$/s & $3\times 10^{-16}$ &  ${}$\cite{Tredenick2024multilayer}  \\ 
		$F$& Faraday constant& sA/mol & 96485.33 &   ${}$\cite{NewellFaraday2019international}  \\ 
		$I$& Current density & A/m$^2$ & $I=I_{\scriptscriptstyle \text{c}}/A$ &    \\ 
		$I_{\scriptscriptstyle \text{c}}$& Applied current for 1C & mA & $5.76$  & Experiment   \\ 
		$k_{\scriptscriptstyle \text{p1}}$& Reaction rate constant & m$^{2.5}$s$^{-1}$mol$^{-0.5}$ & $1.0\times 10^{-10}$ & ${}$\cite{Tredenick2024multilayer}   \\ 
		$k_{\scriptscriptstyle \text{p2}}$& &m$^{2.5}$s$^{-1}$mol$^{-0.5}$ & $8.0\times 10^{-13}$ & ${}$\cite{Tredenick2024multilayer}   \\ 
		$L_{\scriptscriptstyle \text{p1}}$&Thickness & $\mu$m & $44.0$ & Experiment   \\  		
		$L_{\scriptscriptstyle \text{p2}}$& & $\mu$m & $44.0$ & Experiment   \\ 
		$L_{\scriptscriptstyle \text{s}}$& & $\mu$m & $16.0$ & Experiment   \\ 		
		$R$       &  {  Universal gas constant}       & J/K/mol &	{  8.314}         &  \\ 	
		$R_{\scriptscriptstyle \text{c}}$&Contact resistance & $\Omega$ m$^2$ &  $1.5\times 10^{-3}$ & ${}$\cite{Tredenick2024multilayer} \\ 
		$R_{\scriptscriptstyle \text{p1}}$ &Radius of particle & $\mu$m & $4.94$ &   ${}$\cite{Tredenick2025CTvsDFN}   \\ 
		$R_{\scriptscriptstyle \text{p2}}$ &  & $\mu$m & $0.43$ & ${}$\cite{Tredenick2025CTvsDFN}   \\ 
		$T$	&	Constant absolute reference temperature& K & 293.15 ($20\,^{\circ}\mathrm{C}$)& Experiment   \\ 
		$U_0$	&Initial voltage, charge and discharge	 &   V&  3, 4.19& Experiment   \\ 
		$t^+$	&	Transference number of the electrolyte & - & 0.37 &  ${}$\cite{Valoen2005}   \\ 
		$\varepsilon_{\scriptscriptstyle \text{CBD,p1}}$&Carbon binder domain volume fraction& - & 0.11 & ${}$\cite{Tredenick2024multilayer}    \\ 
		$\varepsilon_{\scriptscriptstyle \text{CBD,p2}}$& & - & 0.11 & ${}$\cite{Tredenick2024multilayer}    \\  		
		$\varepsilon_{\scriptscriptstyle \text{e,p1}}$&Electrolyte volume fraction & - & 0.31 & Experiment and  ${}$\cite{Tredenick2025CTvsDFN}    \\  
		$\varepsilon_{\scriptscriptstyle \text{e,p2}}$& & - &0.263 & Experiment and  ${}$\cite{Tredenick2025CTvsDFN}    \\ 
		$\varepsilon_{\scriptscriptstyle \text{e,s}}$& & - & 0.45 & Experiment   \\  		
		$\sigma_{\scriptscriptstyle \text{s,p1}}$&Electronic conductivity in solid & S/m & 5.0 & Experiment   \\ 
		$\sigma_{\scriptscriptstyle \text{s,p2}}$& & S/m & 5.0 & Experiment   \\ 
		\hline 
	\end{tabular}
\end{table}

\section*{Results and Discussion}

\subsection*{Default Benchmark Comparison - 3C Charge  }

\begin{figure} [h!]
	\centering
	\includegraphics[width=0.55\textheight,keepaspectratio]{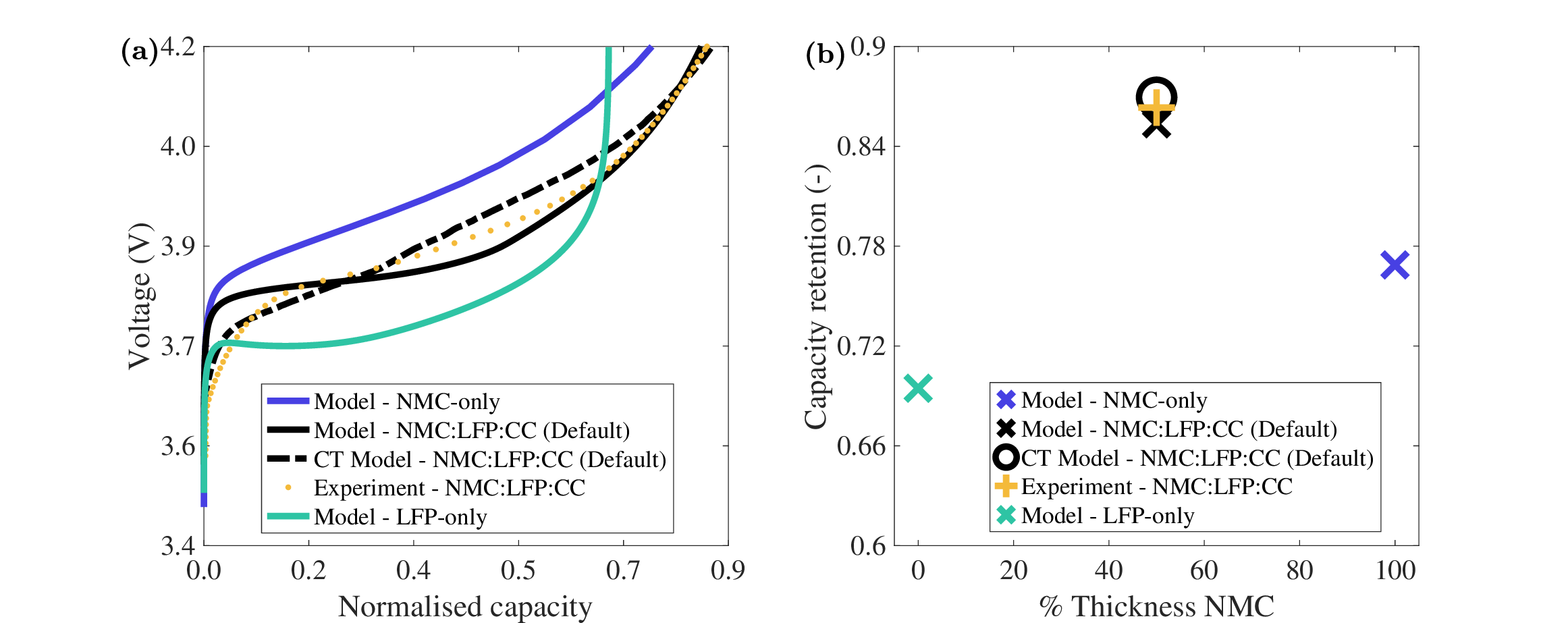}
	\caption{3C charge results with a benchmark comparison including identical specific capacity at 0.05C ($\sim$3.74~mAh/cm$^2$). Conventional half cell electrodes of NMC-only and LFP-only are compared to the bilayer, for voltage as a function of normalised capacity retention (capacity/specific capacity at 0.05C) (a). The capacity retention at 4.2~V is re-plotted in (b). Parameters are shown in Table~\ref{VariablesModel} and \ref{BilNMCLFP} and the experimental data is included from previous work\cite{Tredenick2024multilayer}. The CT model results for the bilayer are also shown, fit to the experimental data and parameters are shown in Table \ref{optimaltable4}.		}
	\label{fig:newcapV3C}
\end{figure}

We begin with a benchmark comparison of the default NMC:LFP:CC bilayer against conventional NMC-only and LFP-only electrodes that have identical specific capacity at 0.05C. The default bilayer is identical to the previous work \cite{Tredenick2024multilayer}, while the conventional electrodes are slightly thicker to closely match the specific capacities at 0.05C ($\sim$3.74~mAh/cm$^2$). The voltage profile during a 3C charge is shown in Fig.~\ref{fig:newcapV3C} as a function of normalised capacity retention (defined as achieved capacity at 4.2~V divided by specific capacity at 0.05C). The specific and normalised capacities are shown in Table~\ref{3CResults}. The bilayer achieves much higher capacity and capacity retention, with LFP-only electrode producing the lowest capacity. The difference between the bilayer and conventional cells capacity retention is 8.5\% and 16\%, for NMC-only and LFP-only, respectively. Fig.~\ref{fig:newcapV3C} also includes the CT model voltage profile, which is validated against the experimental data\cite{Tredenick2024multilayer}. The CT model profile compares well to the experimental data along with the M-DFN modelling results throughout the profile (a), along with at 4.2~V (b).

\begin{table} [h!]
	\footnotesize
	\centering
	\caption{Model parameters for Fig.~\ref{fig:newcapV3C}, for the NMC:LFP:CC bilayer M-DFN and CT, NMC-only, and LFP-only electrode, respectively. }	 	\label{BilNMCLFP}	  
	\begin{tabular}{llllll}
		\hline
		Case        &      & NMC:LFP:CC Bilayer M-DFN & NMC:LFP:CC CT     & NMC-only     & LFP-only      \\  	\hline
		NMC thickness  & $\mu$m       &44 &44 &72 &- \\
		LFP thickness  & $\mu$m          & 44& 44  & -&113 \\ 
		Specific areal capacity 0.05C & mAh/cm$^2$ & 3.741  &3.712 &3.739 &3.740   \\
		Applied current at 1C & mA & 5.762 &2.67$\times 10^{-6}$&5.759   & 5.759  \\	
		Cathode mass & g & 0.035  & 1.5183$\times 10^{-8}$   &0.033   &0.038 \\ \hline
	\end{tabular}
\end{table}

\begin{table}[h!]
	\footnotesize
\centering
	\caption{The specific capacities at 0.05C and achieved capacities at 4.2~V during a 3C charge, for the NMC:LFP:CC bilayer M-DFN and CT,  NMC-only, and LFP-only electrode, respectively.  }
	\label{3CResults}	
	\begin{tabular}{llllll}
		\hline
		        & Units         & NMC:LFP:CC M-DFN & NMC:LFP:CC CT   & NMC-only             & LFP-only\\	\hline 	
		Electrode thickness       & $\mu$m & 88 & 88   & 72    &  113   \\
		Specific capacity at 0.05C     & mAh$/$cm$^2$ & 3.74  & 3.71&3.74 &3.74   \\ \hline 
		Achieved capacity at 4.2~V and 3C         & mAh$/$cm$^2$&    3.19   &3.25 &    2.87     &   2.60  \\	
Capacity retention at 4.2~V and 3C  	& -     &     85.4\%    &     86.9\%    		   &    76.9\%        &   69.5\% \\
	Change in capacity compared to M-DFN bilayer	& -  & -     &     -  		   &  $\downarrow$8.5\%        &   $\downarrow$15.9\%      \\ \hline 
	\end{tabular}
\end{table}

 \subsection*{Benchmark Sensitivity Analysis }
 
\begin{figure} [h!]
 	\centering
 	\includegraphics[width=0.6\textheight,keepaspectratio]{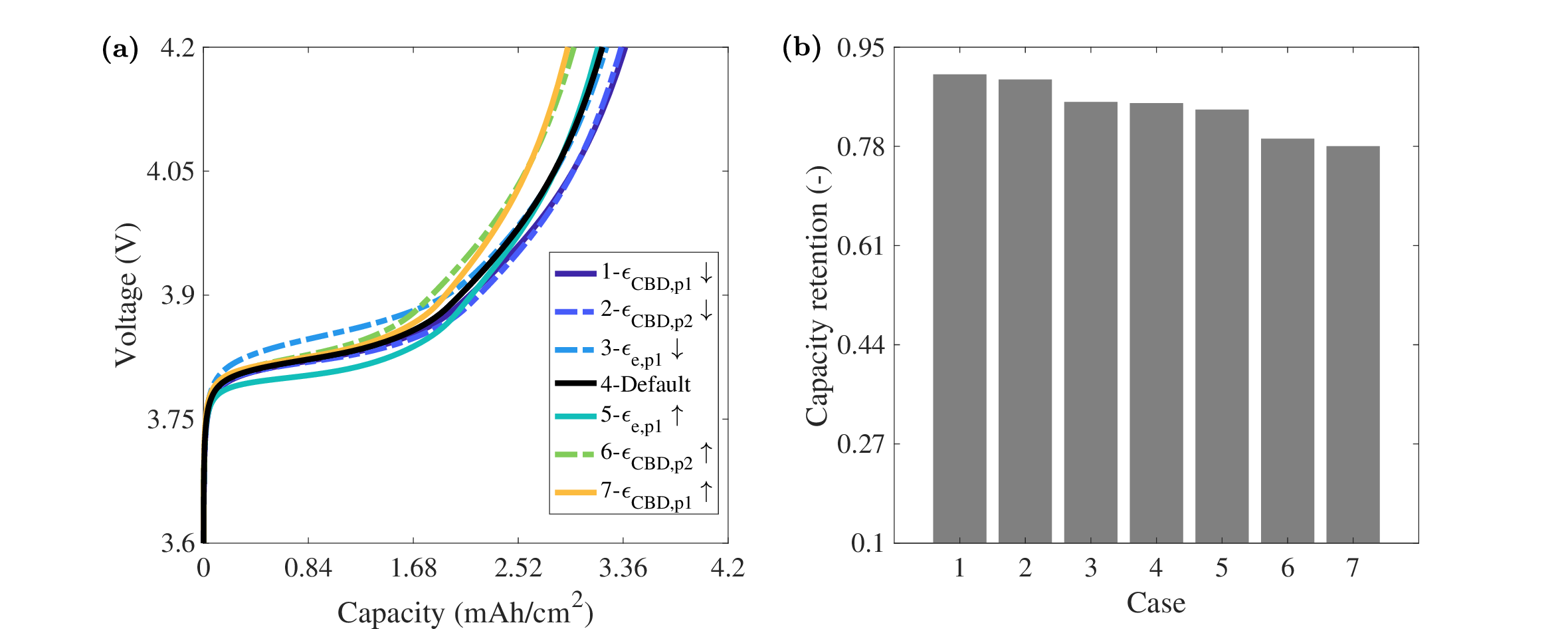}
 	\caption{3C charge sensitivity analysis or benchmark comparison with parameters described in Table \ref{TestsSA3C_2} and the default case 4 shown in black, where p1 is NMC and p2 is LFP.}
 	\label{fig:SA3C_thinner2}
 \end{figure}

 The 3C charge sensitivity analysis for a benchmark comparison with equivalent specific capacities at 0.05C, is shown in Fig.~\ref{fig:SA3C_thinner2} with parameters described in Table \ref{TestsSA3C_2}, including the default case 4. The only parameter that leads to a significant change in the capacity is the volume fraction of the carbon binder domain (CBD), $\varepsilon_{\scriptscriptstyle \text{CBD}}$, as this impacts the area for reactions to occur on the active particle surface. The following parameters were investigated, but not included in Fig.~\ref{fig:SA3C_thinner2}, as no increase in capacity was found, though a decrease was often observed: $\varepsilon_{\scriptscriptstyle \text{e,p2}}, b_{\scriptscriptstyle \text{p1}}, b_{\scriptscriptstyle \text{p2}}, k_{\scriptscriptstyle \text{p1}}$ and $k_{\scriptscriptstyle \text{p2}}$, for the LFP porosity, tortuosity and reaction rate.

 \begin{table} [h!]
 	\scriptsize
 	\centering
 	\caption[Tests]{Sensitivity analysis parameters for Fig.~\ref{fig:SA3C_thinner2}. The parameters have been carefully chosen to obtain a consistent C rate during 3C charging and specific areal capacities at 0.05C ($\sim$3.74~mAh$/$cm$^2$). Case 4 is the default and for cases 3 and 5, the sub-layer thicknesses are changed, along with the porosity, to maintain a consistent specific capacity at 0.05C. The shaded cells are the sensitivity variable and if blank it is equal to the default value of case 4. }	 	\label{TestsSA3C_2}	 
 	\begin{tabular}{lllllllll}
 		\hline
 				Parameter    &Units & Case - 1 &2    & 3    & 4 -Default      & 5 & 6&7     \\	  	 \hline
 		NMC thickness  & $\mu$m & &  &  41 & 44    &53     &     &   \\
 		NMC porosity      & -   & &   & \cellcolor[HTML]{E7E6E6}0.25 & 0.31   &  \cellcolor[HTML]{E7E6E6}0.45    &  &      \\
 		NMC CBD porosity  & - & \cellcolor[HTML]{E7E6E6}0.05  &&   & 0.11&& &\cellcolor[HTML]{E7E6E6}0.2      \\ 	\hline
 		LFP thickness   & $\mu$m   & & & 44 & 44     &  47     &   &     \\
 		LFP  porosity      & -   & &&     & 0.2625 & &  &      \\
 		LFP CBD porosity  & - &  &\cellcolor[HTML]{E7E6E6}0.05 &&0.11   &  &\cellcolor[HTML]{E7E6E6}0.2 &      \\ 	\hline
 		Specific areal capacity at 0.05C & mAh/cm$^2$  & &  & 3.77&  3.74  & 3.75   &   &     \\
 		Applied current at 1C &mA       & &         & 5.81& 5.76   &  5.77  &  &   \\	\hline
 	\end{tabular}
 \end{table}

  \begin{figure} [h!]
 	\centering
 	\includegraphics[width=0.6\textheight,keepaspectratio]{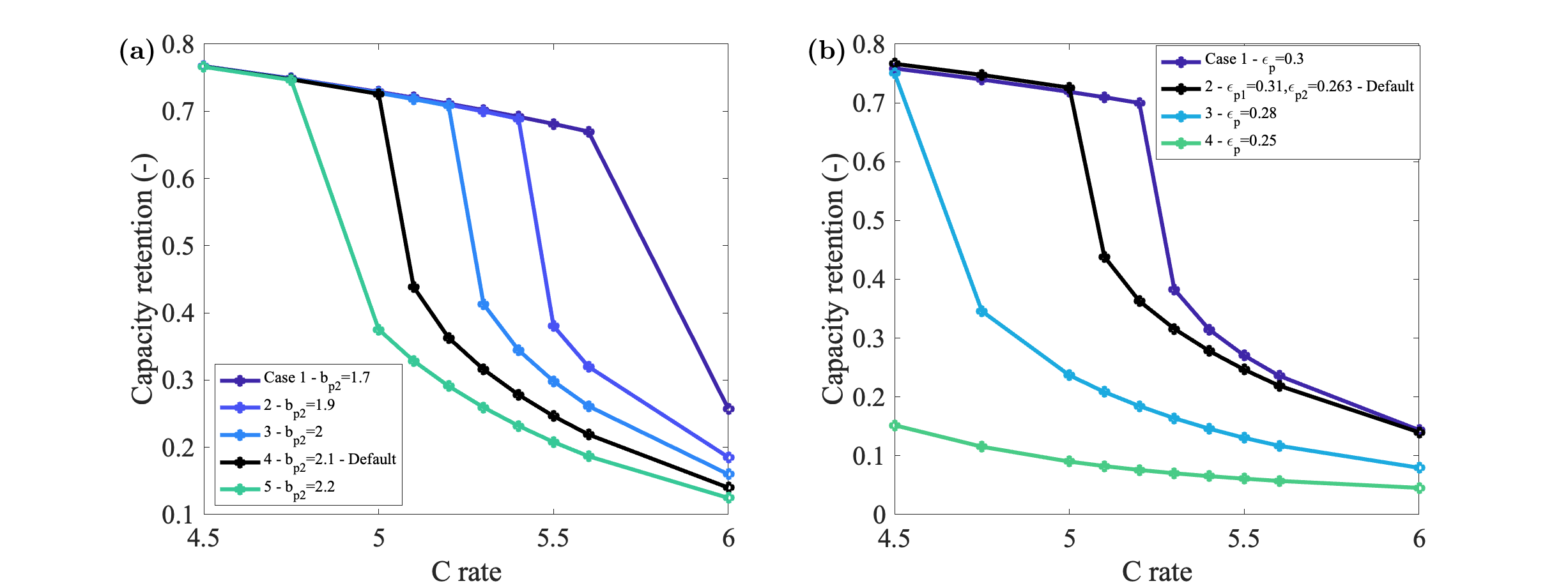}
 	\caption{ The effect of the Bruggeman tortuosity factor ($b_{\scriptscriptstyle \text{p2}}$) of the LFP sub-layer (a) and porosity ($\varepsilon_{\scriptscriptstyle \text{p}}$) (b) on C rate charge performance in terms of capacity retention at 4.2~V. (b) is not a benchmark comparison and the simulation parameters of case 1 are 5.7~mAh and 3.7~mAh/cm$^2$; case 3 are 5.9~mAh and 3.8~mAh/cm$^2$; and case 4 are 6.1~mAh and 4~mAh/cm$^2$. In (b) the two sub-layers of the bilayer have the same porosity except in case 2 for the default.}
 	\label{fig:SAallC}
 \end{figure}

Certain parameters have little effect on the capacity at a given C rate, but can play significant roles in shifting the `capacity cliff', defined as the point where the capacity decreases significantly and is the result of electrolyte depletion and under utilisation of the active materials. Fig.~\ref{fig:SAallC} shows the capacity retention with increasing C rate for a range of tortuosity values for the LFP sub-layer (a), along with a range of porosities (b). As shown in Fig.~\ref{fig:SAallC} (a), when the C rate is increased close to the capacity cliff, a less tortuous diffusion pathway (lower value of the Bruggeman tortuosity factor, $b_p$) can increase capacity, as electrolyte depletion is reduced. At 5.6 C the capacity is significantly higher at  $b_{\scriptscriptstyle \text{p2}}=1.7$ compared to 1.9 as the cliff has shifted. However, if considering only integer C rates (5, 6, 7 C) as is common in practice, then there is not a significant difference. A similar trend is seen in Fig.~\ref{fig:SAallC} (b), where a more porous diffusion pathway (larger porosity value, $\varepsilon_e$, and larger pores), increases capacity. Since high electrolyte potential and electrolyte depletion is occurring beyond the capacity cliff, smaller pores and higher energy density (smaller porosity value, $\varepsilon_e$, for example 0.25) reduces capacity, as this causes the battery to reach its cut-off voltage earlier. If the porosity value is too large, for example $\varepsilon_{\scriptscriptstyle \text{p}}=0.35$ (results not shown for brevity), although the capacity increases, the larger pores reduce the overall energy density of the cell as there is less active material, which makes it unsuitable for EV applications. We therefore adopt a balanced approach of using a porosity in both sub-layers, $\varepsilon_{\scriptscriptstyle \text{p}}$, of 0.3. The results in Fig.~\ref{fig:SAallC} agree with those found elsewhere\cite{tjaden2018tortuosity,DANNER2016191,Usseglio2018Tort} and highlight the importance of balancing energy density and C rate performance, as electrode design is often a trade-off between these two competing metrics.

  \subsection*{Candidate Optimal Parameters}
 The optimal parameters detailed in Table \ref{optimaltable2} are implemented in the M-DFN model, with results shown in Fig.~\ref{fig:firstoptimal1}. These results build upon the outcomes in Figs.~\ref{fig:SA3C_thinner2} and \ref{fig:SAallC}, showing that reducing the carbon binder domain (CBD) volume fraction increases capacity, while changes to the tortuosity and porosity move the capacity cliff. Fig.~\ref{fig:firstoptimal1} demonstrates the sensitivity of the capacity during a 3C charge to changes in these parameter values. The new parameters were chosen with the aim to make several small changes, while retaining relatively high energy density, for example, porosities, $\varepsilon_{\scriptscriptstyle \text{e}} $, of 0.35 and larger, were not considered. We will now take this new set of candidate optimal bilayer parameters and perform further analysis to increase capacity.

 \begin{table} [h!]
 	\footnotesize
 	\centering
 	\caption{Comparison of original and candidate optimal design bilayer parameters.}	 	\label{optimaltable2}	  
 	\begin{tabular}{llll}
 		\hline
 		Parameter        &             &   Default  & Optimal   \\  	\hline
 		NMC electrolyte volume fraction  &   $\varepsilon_{\scriptscriptstyle \text{e,p1}} $     & 0.31 & 0.3\\
 		LFP electrolyte volume fraction  &   $\varepsilon_{\scriptscriptstyle \text{e,p2}} $     &0.263  & 0.3\\
 		NMC  carbon binder domain volume fraction &   $\varepsilon_{\scriptscriptstyle \text{CBD,p1}}$     &  0.11& 0.04\\
 		LFP  carbon binder domain volume fraction  &   $\varepsilon_{\scriptscriptstyle \text{CBD,p2}} $     &0.11  & 0.07\\
 		LFP Bruggeman tortuosity factor	&$b_{\scriptscriptstyle \text{p2}}$& 2.1& 1.8  \\ 
 		\hline
 	\end{tabular}
 \end{table}

 \subsection*{Electrode Thickness - 3C Charge }

\begin{figure} [h!]
 	\centering
 	\includegraphics[width=0.6\textheight,keepaspectratio]{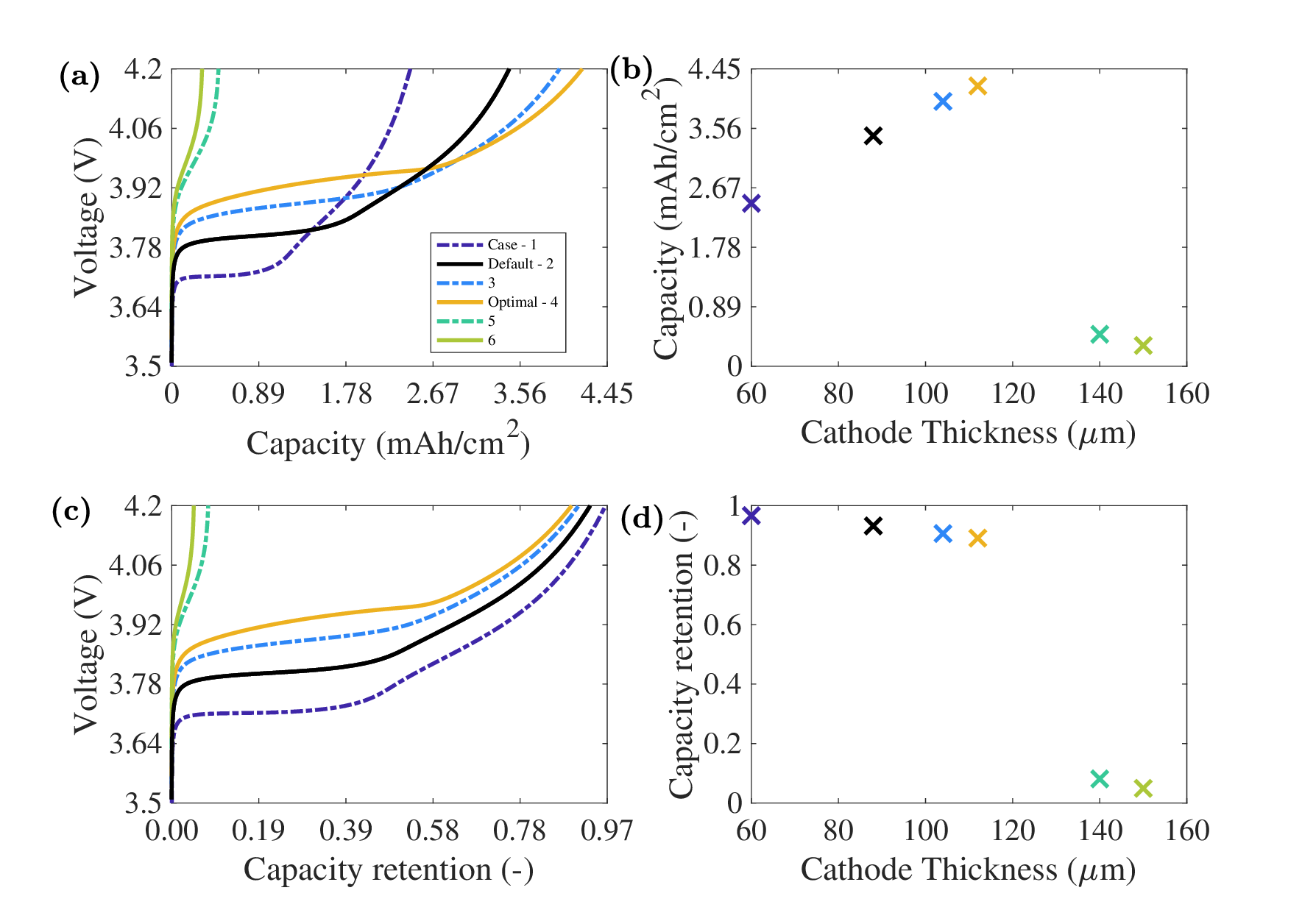}
 	\caption{3C charge with a range of cathode thicknesses. We consider both the capacity (a,b) and capacity retention (c,d) as the specific capacity at 0.05C varies greatly with changes in thickness. The new candidate optimal parameters are shown in Table~\ref{optimaltable2}. Case 2 is the default, case 4 is the new candidate optimal design and thickness changes are shown in Table~\ref{Thickness3C}. Additional thicknesses are shown in Fig.~\ref{fig:thickagain}.}
 	\label{fig:compthtick2_3C}
 \end{figure}

 Fig.~\ref{fig:compthtick2_3C} and Table~\ref{Thickness3C} show the outcome of changing the electrode thickness on the capacity during 3C charging. The voltage profile is shown together with the capacity in  Fig.~\ref{fig:compthtick2_3C} (a), and capacity retention in  Fig.~\ref{fig:compthtick2_3C} (c) (achieved capacity normalised by specific capacity at 0.05C), as the specific capacity at 0.05C varies when varying the thickness. Increasing the thickness increases capacity until an optimal (112~$\mu$m) is reached, then the capacity drops steeply off the capacity cliff for the thickest electrodes. The increase in capacity with the increased electrode thickness is the most substantial gain with any parameter change and the new optimal parameter set helped achieve this by shifting the capacity cliff. In  Fig.~\ref{fig:compthtick2_3C} (d), the capacity retention of the optimal bilayer is slightly lower than the thinner design (case 1) but it is still acceptable (97\% compared to 89\%). 
 
 Fig.~\ref{fig:thickagain} shows additional thickness options, extra to Fig.~\ref{fig:compthtick2_3C}, which exhibit a different profile shape where the voltage increases and then decreases in the middle of the voltage profile, related to a high electrolyte potential and Li ion electrolyte concentration approaching zero (close to $\sim$20~mol/m$^3$). These additional profiles have similar capacity at 4.2~V to the optimal and are likely to be achievable practically, but are not considered further here due to the potential of electrolyte depletion. Fig.~\ref{fig:thickagainCT} for the CT model also exhibit these additional profile shapes but instead by changing the LFP volume fraction. There are several ways to achieve these profile shapes by approaching the capacity cliff including LFP and NMC volume fraction, tortuosity and electrode thickness.

 \begin{table} [h!]
 	\scriptsize
 	\centering
 	\caption{Cathode thickness variations relating to Fig.~\ref{fig:compthtick2_3C}. The default case 2 is shown along with the candidate optimal design of case 4.}	 	\label{Thickness3C}	 
 	\begin{tabular}{llllllll}
 		\hline
 		& Units                    &Case - 1    & Default - 2    & 3 & 4 - Optimal    & 5   & 6          \\ 	\hline
 		NMC thickness       &  $\mu$m & 30     & 44     & 52 &56  &70  & 75      \\
 		LFP thickness       &  $\mu$m   & 30      & 44  & 52  &56 &70 & 75     \\
 		Cathode thickness & $\mu$m  & 60     & 88 &   104 &112 &140   &  150   \\
Applied current at 1C &mA &3.89   & 5.70 & 6.73 &     7.25  &9.07 & 9.71     \\
Specific areal capacity at 0.05C& mAh/cm$^2$  & 2.52    & 3.70 & 4.37 &4.71&5.89    & 6.31   \\	\hline   
 	\end{tabular}
 \end{table}

\subsection*{Compartment Sub-layer Thickness Ratio - 3C Charge }

\begin{figure} [h!]
	\centering
	\includegraphics[width=0.65\textheight,keepaspectratio]{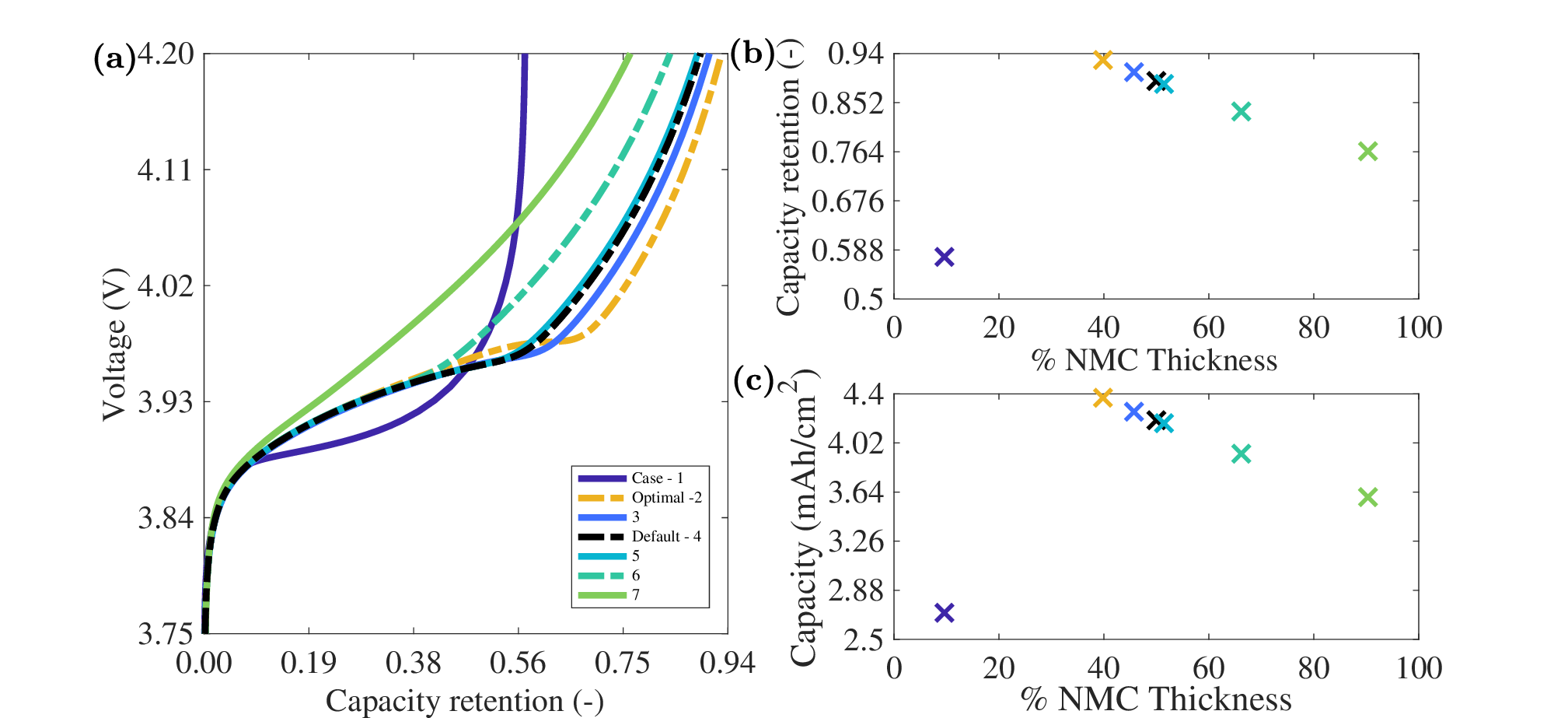} 
	\caption{3C charge with a range of sub-layer thicknesses ratios. This is a benchmark comparison with consistent specific capacities at 0.05C, obtained by varying the total electrode thickness, along with the sub-layer ratio. Case 4 is the rolling default with parameters shown in Table~\ref{CompartThick3C}, and case 2 is the new optimal, including the new parameters from Tables \ref{fig:compthtick2_3C} and \ref{Thickness3C}. }
	\label{fig:compthtick3C}
\end{figure}

Using the optimised electrode thickness from the previous section and the candidate optimal parameters (Tables \ref{fig:compthtick2_3C} and \ref{Thickness3C}), the final step of the M-DFN bilayer design optimisation is to investigate the sub-layer ratio thickness. This is a benchmark comparison and the specific capacity at 0.05C is similar between all cases (4.7~mAh/cm$^2$), as shown in Table~\ref{CompartThick3C}. The results are shown in Fig.~\ref{fig:compthtick3C}, including in (a) the voltage and capacity retention profile, along with the achieved capacity at the cut-off voltage of 4.2~V for both capacity retention, and capacity in (b) and (c). A new optimal is found close to the original 50\% NMC thickness, with a 40\% NMC thickness composed of 47~$\mu$m and 71~$\mu$m thick NMC and LFP sub-layers, respectively. The new optimal produced both higher capacity and capacity retention. 

After tuning the additional parameters of sub-layer thickness ratio and electrode thickness, the new optimal (47~$\mu$m: 71~$\mu$m) had roughly similar NMC sub-layer thickness to the original (44~$\mu$m compared to 47~$\mu$m) but the new candidate optimal parameters (Table \ref{fig:compthtick2_3C}) allowed the LFP sub-layer to increase significantly in thickness from 44~$\mu$m to 71~$\mu$m. This implies that we have generally optimised the LFP sub-layer, while the NMC sub-layer was already close to being optimal in the default 50:50 form.

\begin{table} [h!]
 	\scriptsize
\centering
	\caption{3C charge with a range of sub-layer thicknesses for Fig.~\ref{fig:compthtick3C}. Each case has consistent specific capacities at 0.05C and applied current so is a benchmark comparison.}	 	\label{CompartThick3C}	 
	\begin{tabular}{lllllllll}
		\hline
		&Units  &Case - 1 & 2 - Optimal  & 3& 4 - Default & 5  & 6 & 7     \\ \hline
		NMC thickness &  $\mu$m & 13.50 &47.00  &52.25&  56.00& 57.25 &68.50  &84.00     \\
		LFP thickness &$\mu$m &  127.00&71.00  &62.00  &56.00  &54.00  & 35.00 & 9.00  \\
		NMC thickness & \% & 9.6 &39.8  & 45.7 &50.0  &51.5  & 66.2 &90.3   \\
		Applied current at 1C &mAh  & 7.24 & 7.25 & 7.24 &7.25  & 7.26 &7.25  & 7.25  \\
		Specific areal capacity at 0.05C & mAh/cm$^2$  & 4.70 & 4.71 &  4.70& 4.71 & 4.71 & 4.71 &4.71   \\
		\hline
	\end{tabular}
\end{table}

\subsection*{ Optimal Bilayer - 3C Charge}

\begin{figure} [h!]
 \centering
	\includegraphics[width=0.33\textheight,keepaspectratio]{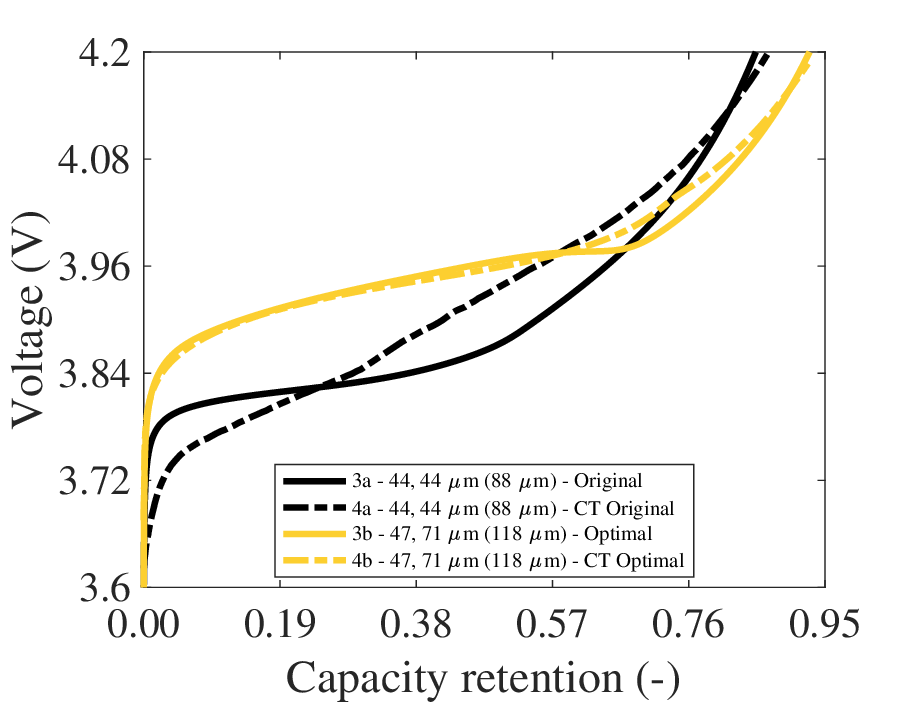}
	\caption{Optimal of case 3b and 4b, compared to the original bilayer of case 3a and 4a, during a 3C charge, for the M-DFN model along with the CT model. The voltage is shown as a function of capacity retention. Cases 3a and 4a have capacity $\sim$3.7~mAh/cm$^2$ and cases 3b and 4b have capacity $\sim$4.7~mAh/cm$^2$. The thicker case 3b includes the candidate optimal design parameters in Tables~\ref{optimaltable2} and \ref{optimaltable}. The capacities at 4.2~V and additional voltage profiles are shown in Figs.~\ref{fig:optimal} (c) and \ref{fig:optimal3}.
	}	
	\label{fig:optimal1}
\end{figure}

\begin{figure} [h!]
 	\centering
	\includegraphics[width=0.74\textheight,keepaspectratio]{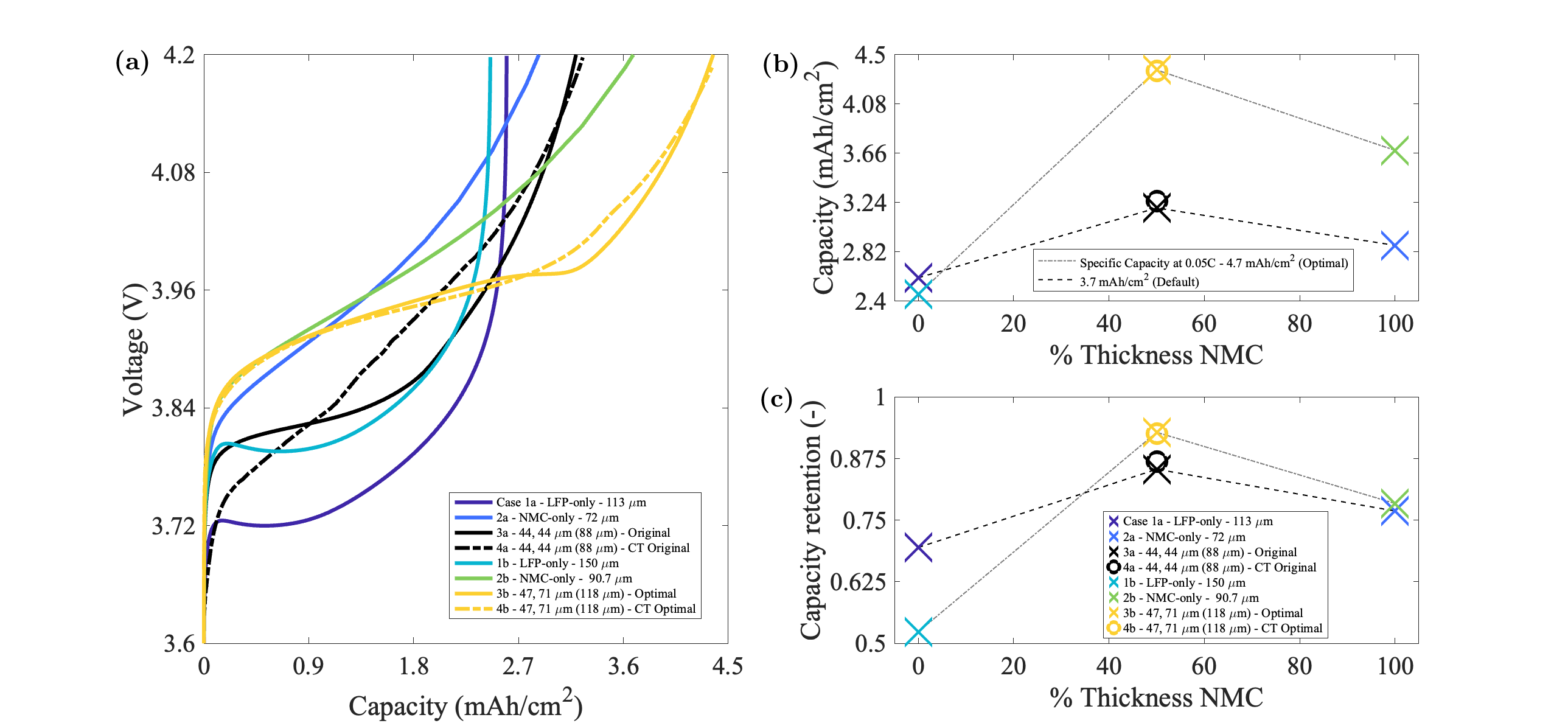}
	\caption{Optimal of case 3b and 4b, compared to the original bilayer of case 3a and 4a, during a 3C charge, for the M-DFN model along with the CT model. The NMC-only and LFP-only cases (1a, 2a, 1b, 2b) are included with equivalent specific capacity at 0.05C. Both the capacity (a) and (b), and normalised capacity (c), are shown. Cases 1a to 4a have capacity $\sim$3.7~mAh/cm$^2$ and cases 1b to 4b have capacity $\sim$4.7~mAh/cm$^2$. The thicker cases 1b to 3b include the candidate optimal design parameters in Tables~\ref{optimaltable2} and \ref{optimaltable}.
	}	
	\label{fig:optimal}
\end{figure}

The optimal bilayer cathode, with parameters shown in Tables~\ref{CompartThick3C}, \ref{optimaltable}, and \ref{optimaltable3}, is compared to the original bilayer design, along with the CT model and conventional NMC-only and LFP-only electrodes, as shown in Figs.~\ref{fig:optimal1}, \ref{fig:optimal} and \ref{fig:optimal3}. Figs.~\ref{fig:optimal1} and \ref{fig:optimal} show voltage as a function of capacity along with capacity retention, and compare equivalent specific capacities at 0.05C of $\sim$3.7~mAh/cm$^2$ for the default and $\sim$4.7~mAh/cm$^2$ for the optimal. Fig.~\ref{fig:optimal1} shows the original and optimal bilayer of the M-DFN model compared to the CT model for capacity retention. The two models compare well in terms of both achieved capacity at 4.2~V and the shape of the voltage profiles, indicating that the CT model has agreed well with optimal found using the M-DFN model. Fig.~\ref{fig:optimal} shows the thicker conventional NMC-only and LFP-only electrodes in cases (1b, 2b) also include the new optimised parameters including porosity, CBD and tortuosity from Table~\ref{optimaltable2}, for a fair comparison to the optimal bilayer and to maintain equivalent capacity at 0.05C. The new optimal bilayer has increased capacity retention with significantly increased capacity in mAh/cm$^2$ compared to the NMC-only and LFP-only electrodes, as shown in Table~\ref{3CResults2}. With small changes to five parameters (Table \ref{optimaltable2}), along with larger changes to the electrode thickness and sub-layer thickness, we were able to significantly increase capacity in mAh/cm$^2$, while still maintaining reasonable and improved capacity retention.

In Fig.~\ref{fig:optimal} and Table \ref{3CResults2}, comparing the conventional NMC-only (case 2b) and LFP-only (case 1b) electrodes capacity to the optimal bilayered electrode (case 3b) with equivalent specific areal capacity, the new optimal bilayer increases capacity retention during a 3C charge by 15\% or 0.7~mAh/cm$^2$, and 41\% or 1.9~mAh/cm$^2$, respectively. The new optimal bilayer charges between 0-90\% SOC in 18.6 minutes, achieving 4.4~mAh/cm$^2$. The NMC-only electrode charges between 0-78\% SOC in 15.7 minutes but only achieves 3.7~mAh/cm$^2$. The LFP-only electrode charges between 0-52\% SOC in 10.5 minutes but only achieves 2.5~mAh/cm$^2$. Compared to the original bilayer (case 3a), the optimal bilayer (case 3b) achieves a slightly higher capacity retention of 7.5\%, along with a 1.2~mAh/cm$^2$ higher capacity and an increase in the electrode thickness from 88~$\mu$m to 118~$\mu$m.

\begin{table} [h!]
	\footnotesize
	\centering
	\caption{Parameters for Fig.~\ref{fig:optimal}.}	 	\label{optimaltable}	  
	\begin{tabular}{lllllll}
		\hline
Specific capacity 0.05C	&	        &    &1a - LFP-only    & 2a - NMC-only     &  3a - M-DFN  Default    &  4a - CT Default   \\  	\hline
3.7~mAh/cm$^2$		&	NMC thickness  & $\mu$m   &-&72    &44 &44\\
&	LFP thickness  & $\mu$m      &113     & - & 44& 44\\ 
	&	Applied current at 1C & mA &5.76 &5.76 & 5.76 & 2.67$\times 10^{-6}$ \\	 \hline 
	&		       &     &   1b - LFP-only  &2b - NMC-only & 3b - M-DFN Optimal & 4b - CT Optimal    \\  	\hline
		4.7~mAh/cm$^2$			&	NMC thickness  & $\mu$m   &-&89.2 &47 &47\\
	&		LFP thickness  & $\mu$m    &149.5&-& 71& 71\\ 
		&		Applied current at 1C & mA  &7.23 & 7.24 &7.25 &3.21$\times 10^{-6}$ \\	\hline 
	\end{tabular}
\end{table}

\begin{table}[h!]
	\scriptsize  
 	\centering
	\caption{The specific and achieved capacities at the end of charge at 4.2~V for Fig.~\ref{fig:optimal}.}
	\label{3CResults2}	
	\begin{tabular}{ | m{2cm} | m{5.8cm}| m{1cm} | m{1.6cm} | m{1.8cm}| m{1.5cm} |  m{1.5cm}  |} 
		\hline
Specific capacity	&			      &   &   1a - LFP-only   & 2a - NMC-only     &  3a - M-DFN Default   &  4a - CT Default  \\  	\hline
	3.74~mAh/cm$^2$	&	Achieved capacity at 4.2~V and 3C         & mAh$/$cm$^2$&2.60 &   2.87  &3.19  &3.25 \\	
	&	Normalised capacity at 4.2~V and 3C  	& -  &  69.5\%  &76.9\%    & 85.4\%  & 86.9\% \\	
	&	Change in capacity compared to M-DFN bilayer	&   mAh$/$cm$^2$&  $\downarrow$0.6 &    $\downarrow$0.3&-    &- \\ 
	&	Change in capacity compared to M-DFN bilayer	&-&$\downarrow$15.9\%  &  $\downarrow$8.5\% &-&- \\ 	\hline 
	 	&			      &   &1b - LFP-only  & 2b - NMC-only & 3b - M-DFN Optimal & 4b - CT Optimal    \\  	\hline
	4.7~mAh/cm$^2$	&	Achieved capacity at 4.2~V and 3C         & mAh$/$cm$^2$& 2.46  & 3.68 & 4.37  & 4.36  \\	
	&	Normalised capacity at 4.2~V and 3C  	& -   &52.3\%&   78.4\%&92.9\%  & 92.6\% \\	
	&	Change in capacity compared to M-DFN bilayer	&   mAh$/$cm$^2$&     $\downarrow$1.9 &  $\downarrow$0.7	&- &-  \\ 
	&	Change in capacity compared to M-DFN bilayer	&- &     $\downarrow$ 40.6\% &   $\downarrow$14.5\%	&- &-  \\ \hline
	\end{tabular}
\end{table}

 \begin{table} [h!]
 	\footnotesize
 	\centering
 	\caption{Comparison of original and final optimal bilayer M-DFN design parameters.}	 	\label{optimaltable3}	  
 	\begin{tabular}{lllll}
 		\hline
 		Parameter        &    &      Unit   &   Default  & Optimal   \\  	\hline
 		 		NMC sublayer thickness     &  &  $\mu$m  & 44 &47    \\
 		 		LFP sublayer thickness     &  &  $\mu$m  & 44 &71   \\
 		 		Cathode electrode thickness     &  &  $\mu$m  & 88 &118  \\
 		 	NMC thickness      &  &  \% &50.0 &39.8  \\ 
 		 Applied current at 1C     &  &  mAh  &5.76 &7.25  \\
 		 Specific areal capacity at 0.05C    &  &  mAh/cm$^2$  &3.74 &4.71    \\ \hline
 		NMC electrolyte volume fraction  &   $\varepsilon_{\scriptscriptstyle \text{e,p1}} $     &- & 0.31 & 0.30\\
 		LFP electrolyte volume fraction  &   $\varepsilon_{\scriptscriptstyle \text{e,p2}} $     &- &0.26  & 0.30\\
 		NMC  carbon binder domain volume fraction &   $\varepsilon_{\scriptscriptstyle \text{CBD,p1}}$     &   -&0.11& 0.04\\
 		LFP  carbon binder domain volume fraction  &   $\varepsilon_{\scriptscriptstyle \text{CBD,p2}} $     & -&0.11  & 0.07\\
 		LFP Bruggeman tortuosity factor	&$b_{\scriptscriptstyle \text{p2}}$&-& 2.1& 1.8  \\ 
 		\hline
 	\end{tabular}
 \end{table}

\subsection*{ Optimal Bilayer -  Modelling Results}

\begin{figure*}[h!]
	\centering
	\begin{subfigure}[t]{0.5\textwidth}
		\centering
		\includegraphics[height=2.7in]{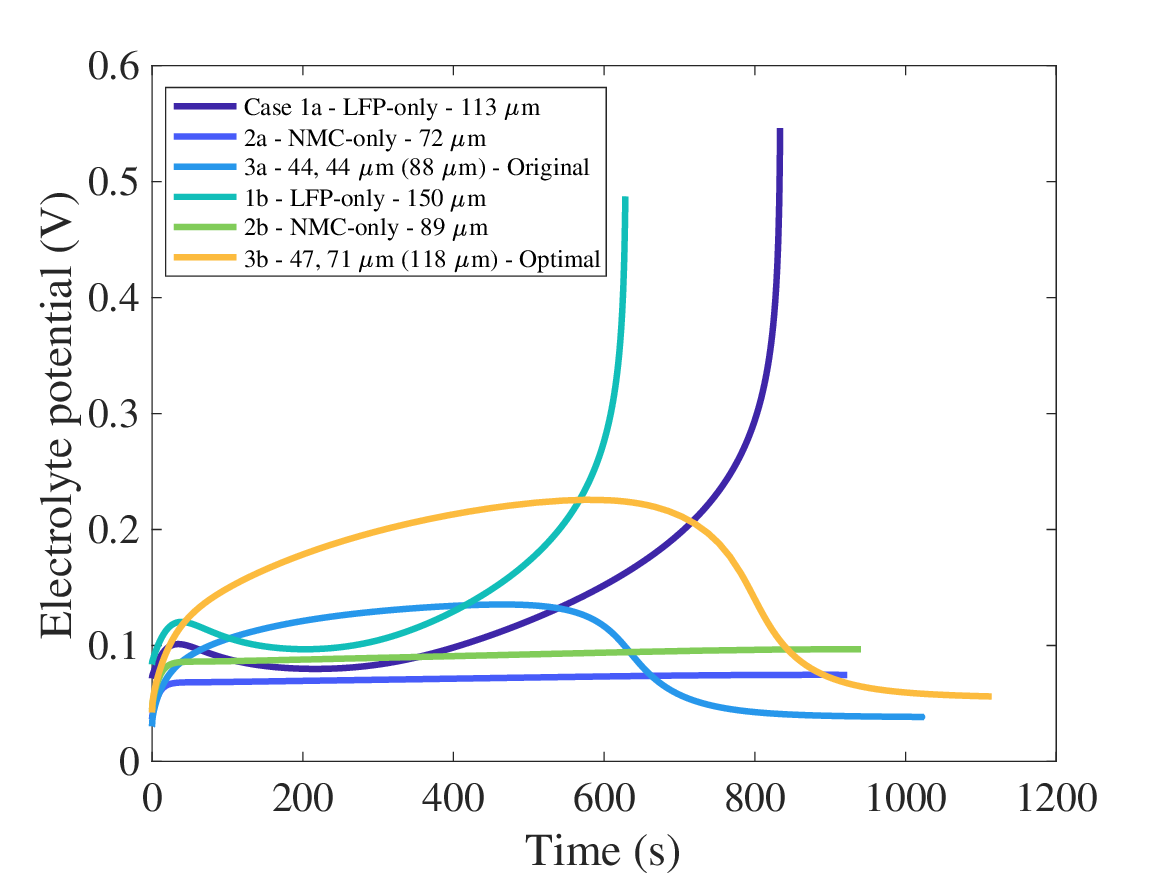}
		\caption{Electrolyte potential.}
	\end{subfigure}%
	~ 
	\begin{subfigure}[t]{0.5\textwidth}
		\centering
		\includegraphics[height=2.7in]{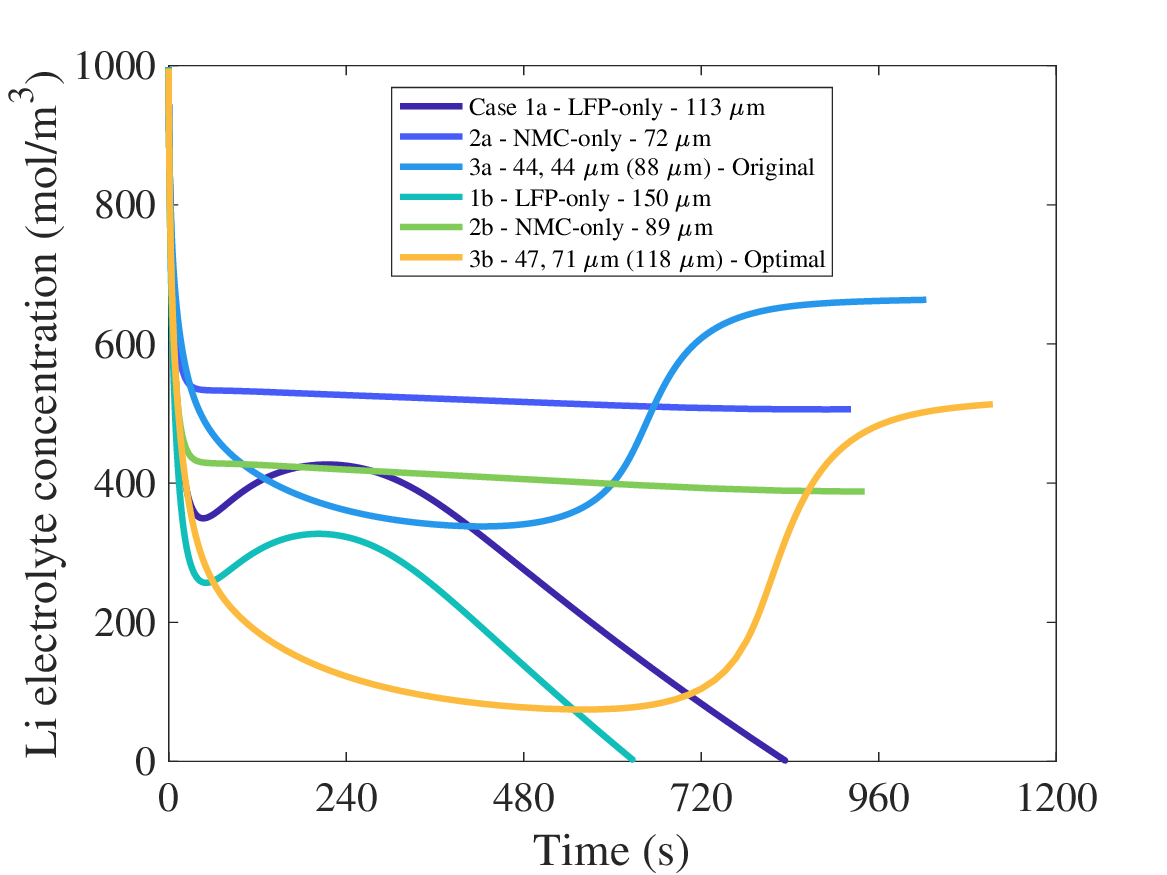}
		\caption{Li ion electrolyte concentration.}
	\end{subfigure}
	\caption{The Li ion electrolyte potential adjacent to the current collector (CC) as a function of time is shown in (a), with the legend showing cases corresponding to Fig.~\ref{fig:optimal}. The Li ion electrolyte concentration adjacent to the separator at the counter electrode (CE) as a function of time is shown in (b). }
	\label{fig:phie} 	\label{fig:ces}
\end{figure*}

\begin{figure} [h!]
	\centering
	\includegraphics[width=0.65\textheight,keepaspectratio]{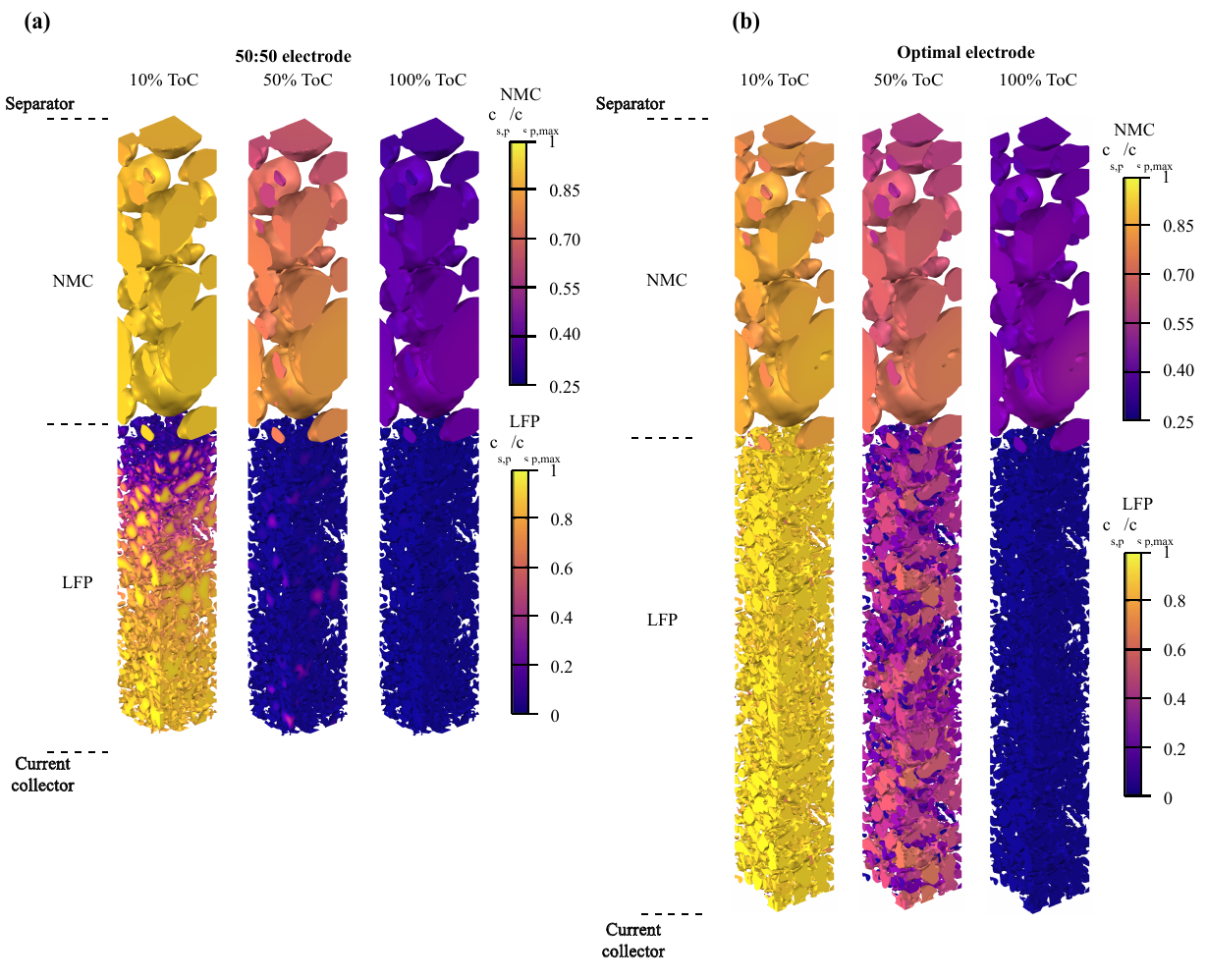} 
	\caption{3D CT model results for normalised Li$^+$ concentration in the active particles ($c_{\scriptscriptstyle \text{s,p}} /c_{\scriptscriptstyle \text{s,p,max}}$) for the default 50:50 bilayer (a) and optimal electrode (b). The Time of Charge (ToC) is shown near the start (10\%), half way (50\%) and at the end of charge (100\%). The NMC particles are shown adjacent to the separator and LFP adjacent to the current collector.}	
	\label{fig:3d_CT}
\end{figure}

 Figs.~\ref{fig:phie} (a) and (b) show the Li ion electrolyte potential adjacent to the current collector and electrolyte concentration adjacent to the separator, corresponding to the cases in Fig.~\ref{fig:optimal}. Upon further inspection of the results for each electrode, including the Li ion concentration and potential in the electrolyte, the two bilayers were observed to be similar, but the new parameter set allowed a thicker electrode, while maintaining the same mechanisms that promote higher capacity, discussed in detail previously\cite{Tredenick2024multilayer}. The potential of the optimal bilayer of case 3b is in between the others so balances the electrolyte response and does not exhibit electrolyte depletion, promoting higher capacity. Conversely, the lower capacity of the LFP-only electrodes of cases 1a and 1b is caused by the following process: \textit{i}) excessive electrolyte potential towards the end of charge (Fig.~\ref{fig:phie} (a)), which causes \textit{ii}) the electrolyte concentration to drop to zero (electrolyte depletion) (Fig.~\ref{fig:phie} (b)), which causes \textit{iii}) the voltage to reach the cut-off voltage prematurely, that in turn causes under utilisation of the active materials and a lower state-of-charge\cite{Tredenick2024multilayer}. The relatively thick LFP-only electrodes at this C rate would be operating outside the manufacturers recommendations due to the reduced capacity, electrolyte depletion and high potential. The NMC-only (case 2b) does not suffer electrolyte depletion problems but is too thick for this C rate to reach a high SOC at 4.2~V, whereas the bilayer benefits from the 0-100\% SOC (see Fig.~\ref{fig:3d_CT}) of the LFP sub-layer\cite{Tredenick2024multilayer}.

\begin{figure} [h!]
	\includegraphics[width=0.67\textheight,keepaspectratio]{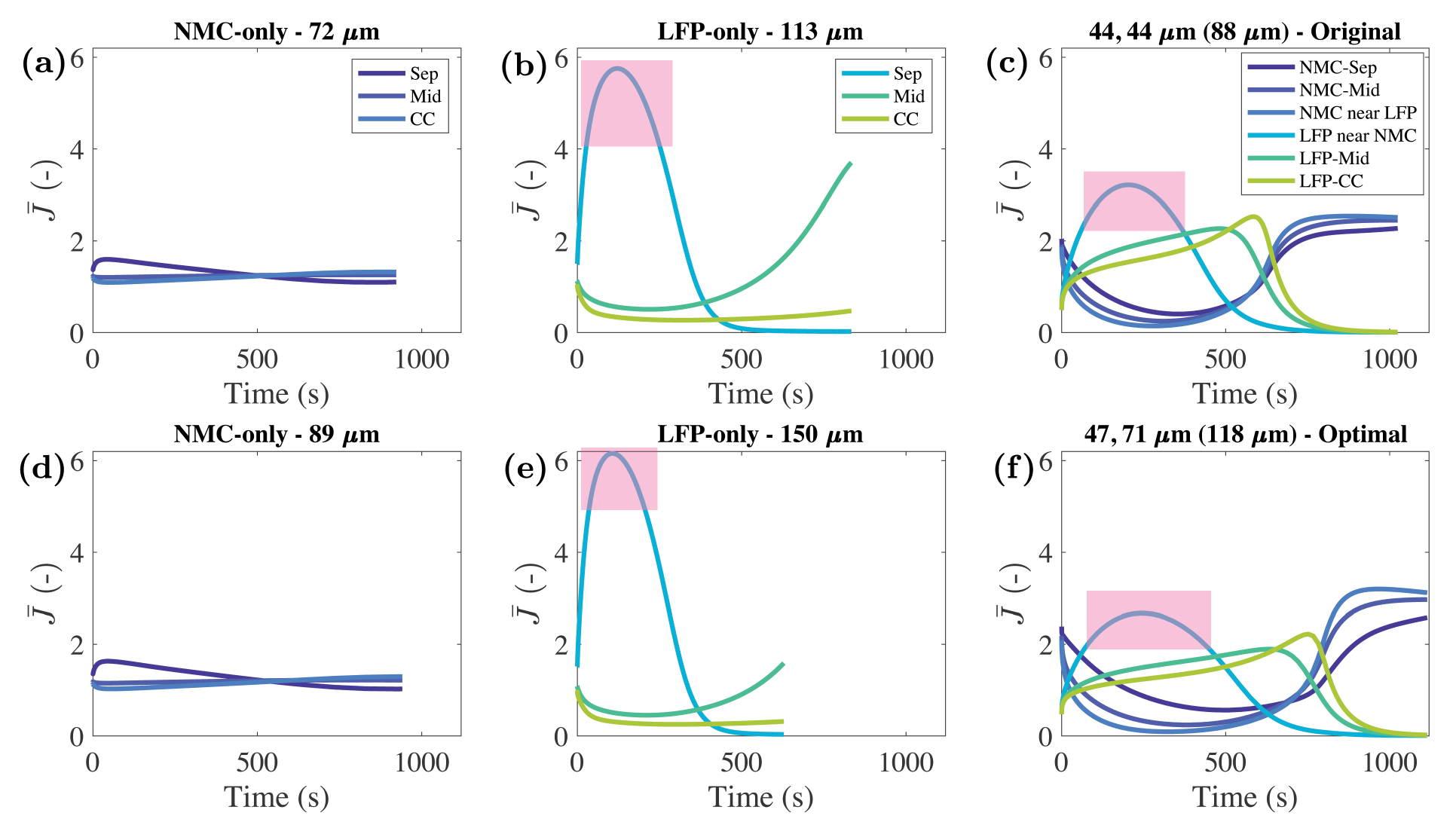} 
	\caption{The normalised Butler-Volmer reaction current density ($\bar{J} = J_{\scriptscriptstyle \text{p,i}}  \  L \   a_{\scriptscriptstyle \text{p,i}}   /  I  $) as a function of time, corresponding to the results in Fig.~\ref{fig:optimal}. The first row (a) to (c) are the original thinner electrodes, with specific capacity at 0.05C of $\sim$3.7~mAh/cm$^2$ and the second row (d) to (f) $\sim$4.7~mAh/cm$^2$. The highlighted regions are the LFP layer peak current densities that are impacted by the increased thickness and bilayer arrangement. The legend shows the current density location including adjacent to the counter-electrode and separator (Sep), in the middle of the bilayer sub-layer or electrode (Mid), and adjacent to the current collector (CC).}
	\label{fig:Jbar}
\end{figure}

Fig.~\ref{fig:Jbar} shows the reaction current density as a function of time, for the cases corresponding to Fig.~\ref{fig:optimal}. The reaction current density is normalised by the applied current, reaction surface area, and electrode thickness, which allows direct comparison across different electrode compositions. The conventional electrodes are shown in (a, b, d, e) and the NMC-only electrodes exhibit small changes in current while the LFP-only has high current next to the separator, and the thicker LFP-only electrode in (e) has the highest current. The LFP-only sub-layer of the bilayer (c, f) compared to LFP-only electrodes (b,e) has lower peak current, as highlighted in pink. Compared to the original in (c), the optimal bilayer in (f) shows that some of the current has shifted from the LFP sub-layer to the higher capacity NMC sub-layer and can charge both layers more homogeneously to promote higher capacity, which could be beneficial in reducing degradation in the LFP sub-layer.

Fig.~\ref{fig:Jbar} shows that, along with the NMC sub-layer thickness not changing significantly in the optimal compared to the default, that we optimised mainly the LFP sub-layer. LFP in general including the overly thick LFP-only electrode (case 1b) is more prone to electrolyte depletion than NMC, especially at high C rates and thicknesses. The high electrolyte potential (Fig .~\ref{fig:phie}) and current (Fig .~\ref{fig:Jbar}) in the LFP-only electrodes could lead to decreased lifespan\cite{xu2019heterogeneous,edge2021lithium} and their use not recommended. The bilayer arrangement enables concentration, potential and current to be in-between the LFP-only and NMC-only electrodes, to promote increased capacity. This shows the benefits of the bilayer approach, where we are able to utilise a higher specific capacity at 0.05C than a conventional LFP-only electrode. 
 
Fig.~\ref{fig:3d_CT} shows the 3D CT model results for normalised Li$^+$ concentration in the active particles ($c_{\scriptscriptstyle \text{s,p}} /c_{\scriptscriptstyle \text{s,p,max}}$) for the default 50:50 bilayer (a) and optimal electrode (b). The NMC particles show uniform concentrations throughout the particle radius and electrode thickness, while the LFP particles show concentration gradients through both the particle radius in (a) and (b) and electrode thickness in (a) due to slower solid state diffusion. The NMC results appear similar between the default and optimal, while the thicker optimal LFP layer takes longer to charge and indicates the optimisation routine resulted in an optimised LFP layer, while NMC was already close to being optimised, as discussed in the previous paragraph.

\subsection*{Optimal Bilayer - C Rate Performance}
\begin{figure} [h!]
	\centering
	\includegraphics[width=0.6\textheight,keepaspectratio]{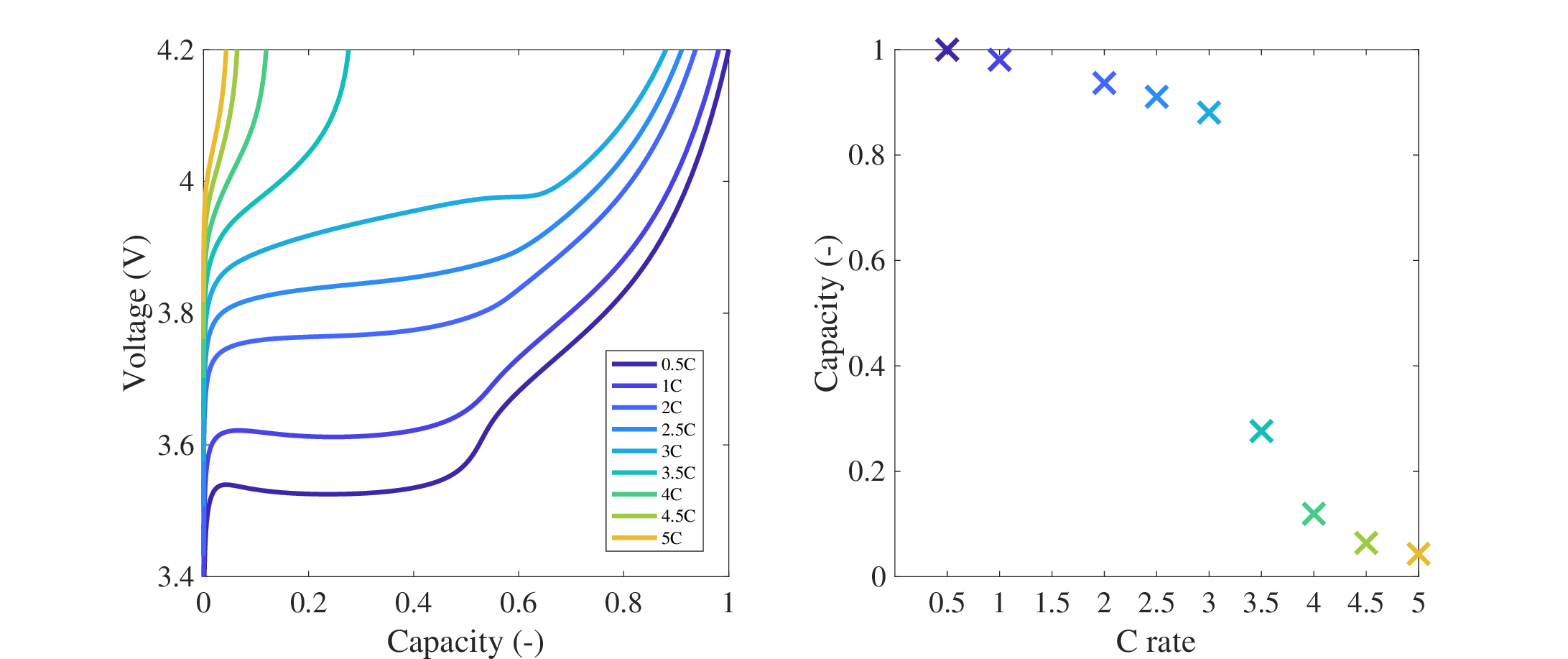}
	\caption{   Charge of the new optimal case including voltage profiles as a function of capacity. A range of C rates is investigated including 0.5, 1, 2, 2.5, 3, 3.5, 4, 4.5 and 5C. The normalised capacity is the achieved capacity at a given C rate divided by the capacity at 0.5C.	}	
	\label{fig:optimal_lotsCrates}
\end{figure}

\begin{figure} [h!]
 	\centering
 	\includegraphics[width=0.35\textheight,keepaspectratio]{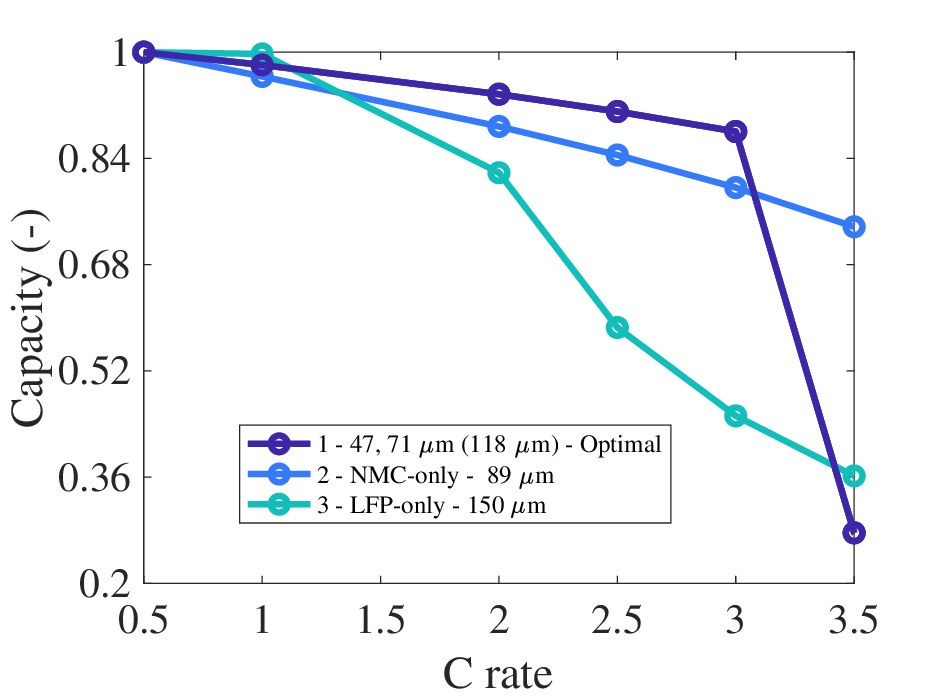}
 	\caption{Normalised capacity for a range of C rates, of the optimal, NMC-only and LFP-only electrodes, at 0.5, 1, 2, 2.5, 3 and 3.5C, corresponding to the results in Fig.~\ref{fig:optimal} of cases 1b, 2b and 3b. The normalised capacity is the capacity at a given C rate divided by the achieved capacity at 0.5C.	}
 	\label{fig:optimalragone2}
 \end{figure}

\begin{figure} [h!]
 	\centering
 	\includegraphics[width=0.38\textheight,keepaspectratio]{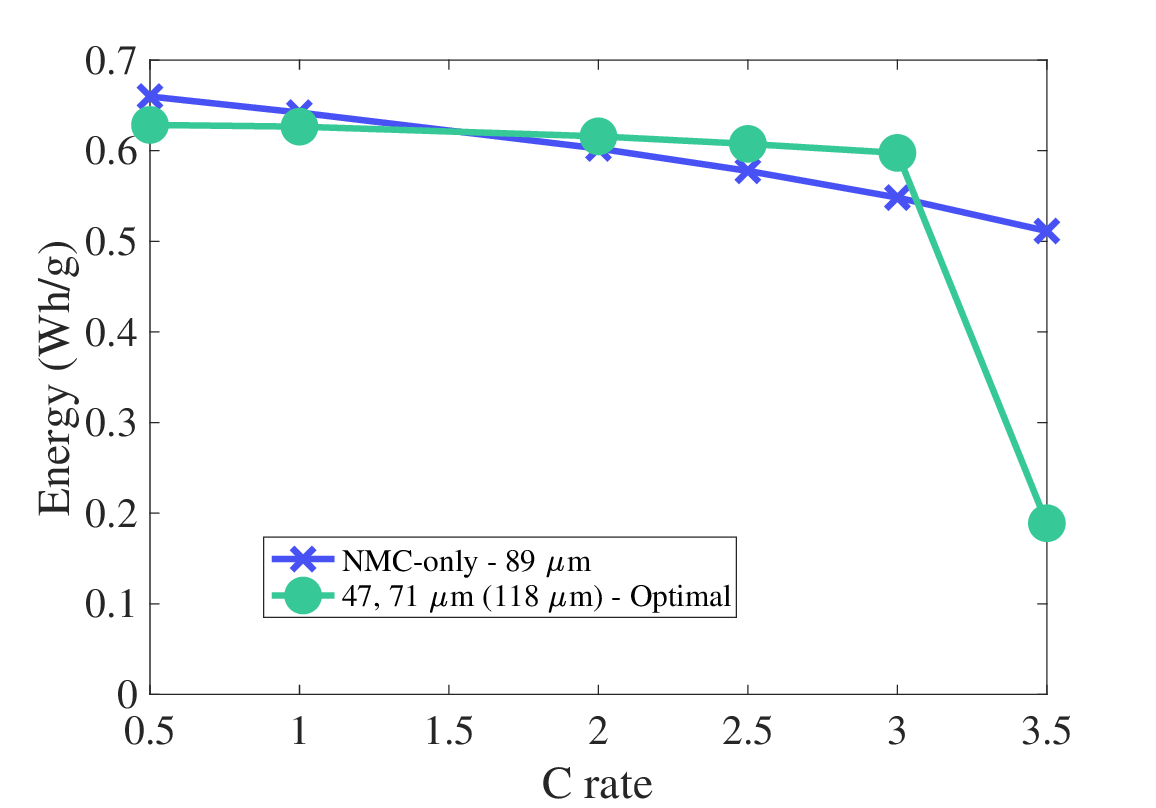}
 	\caption{Mass energy density (Wh/g) for a range of C rates, of the optimal and NMC-only electrodes, at 0.5, 1, 2, 2.5, 3 and 3.5C (LFP-only not shown due to poor performance). The mass is included for each cathode of 44.7~mg for bilayer and 41.4~mg for NMC-only. The energy\cite{beyers2023ragone}, $E$ (Wh/g), is $E= I  / m \int_{t=0}^{t=t_f} V(t) \  dt$, where $V(t)$ is the resulting voltage produced with the model in volts (V) for each case and C rate, $I$ is the applied current (constant in time) in amps (A) for each case and C rate, $t_f$ is the final time in hours (h), and $m$ is the mass of the electrode in grams (g) for each cathode.  	}
 	\label{fig:optimalragonevolume}
 \end{figure}

Fig.~\ref{fig:optimal_lotsCrates} evaluates the optimised bilayer over a range of C rates. The voltage is shown as a function of normalised capacity including C rates up to 5C. As 3C was the target charging rate for the design optimisation, the edge of the capacity cliff was located at 3C by design, and then increasing the C rate higher reduced the capacity. The profile shapes are as expected and explained in detail previously\cite{Tredenick2024multilayer}, including at lower C rates the LFP sub-layer (flat voltage plateau) charges first and separately to the NMC sub-layer (steep voltage profile), which charges in the second half of the profile, due to the different voltage windows of the OCP profiles. At charging rates above 2C, the sub-layers charge together and the voltage profiles appear more smooth. Fig.~\ref{fig:optimalragone2} shows the C rate performance of the optimal bilayer along with the NMC-only and LFP-only electrodes. The excessively thick LFP-only electrode has capacity reduction after 1C, while the optimal and NMC-only electrode have reasonable capacity until 3C. The optimal bilayer has the highest capacity at 2C and 3C and compared to LFP-only, allows the capacity cliff to be shifted to higher C rates. We investigate the energy density dependence with C rate in Fig.~\ref{fig:optimalragonevolume}, for the equivalent NMC-only and optimal bilayer electrodes at C rates 0.5 to 3.5. The NMC-only and bilayer electrodes have similar energy densities at low C rates and at 3C the bilayer has higher energy density as the NMC-only electrode has lower capacity retention (Fig.~\ref{fig:optimal}). 

\begin{figure} [h!]
	\centering
	\includegraphics[width=0.35\textheight,keepaspectratio]{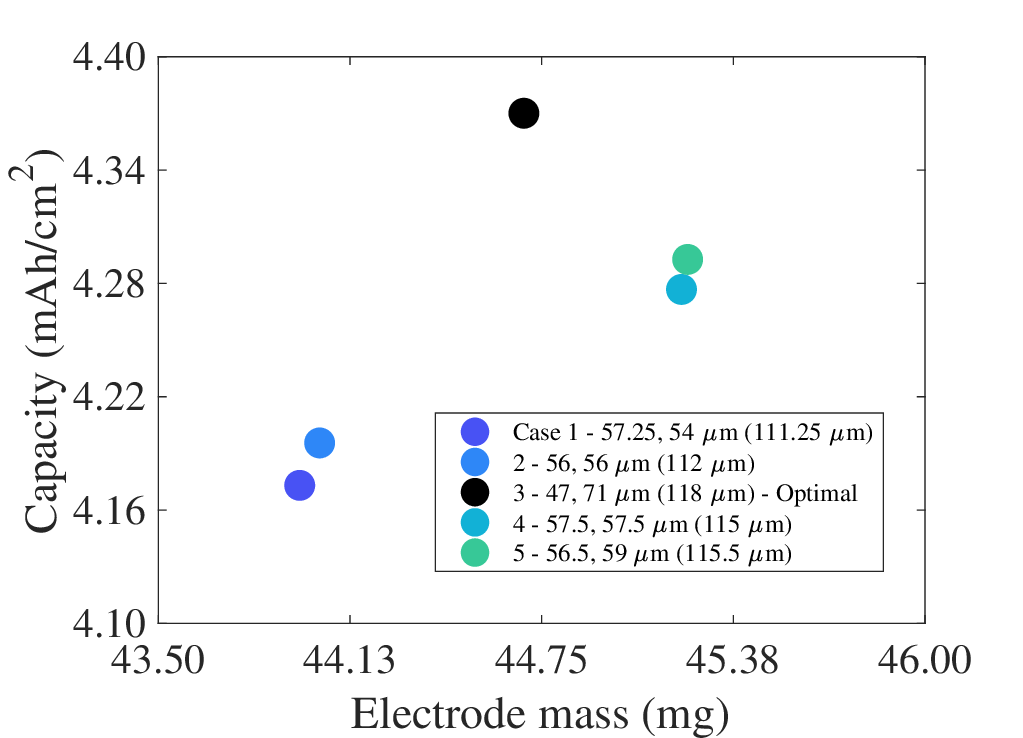}
	\caption{The 3C charge capacity as a function of electrode mass. For cases 1,2,4 and 5 the simulation parameters are 7.26, 7.25, 7.45 and 7.44~mA; 4.71, 4.71, 4.84 and 4.83~mAh/cm$^2$; and 44.0, 44.0, 45.2 and 45.2~mg. }
	\label{fig:withgrams}
\end{figure}

As shown in Fig.~\ref{fig:withgrams}, the active material mass of the electrode is the optimal variable, as apposed to electrode thickness or sub-layer thickness ratio, producing the optimal mass of 44.7~mg. We can arrive at the optimal mass in three ways by varying the porosity (not shown for brevity), electrode thickness or sub-layer ratio thickness.

  \subsection*{Optimal Bilayer - Cycling}
  
  \begin{figure} [h!]
  	\centering
  	\includegraphics[width=0.6\textheight,keepaspectratio] {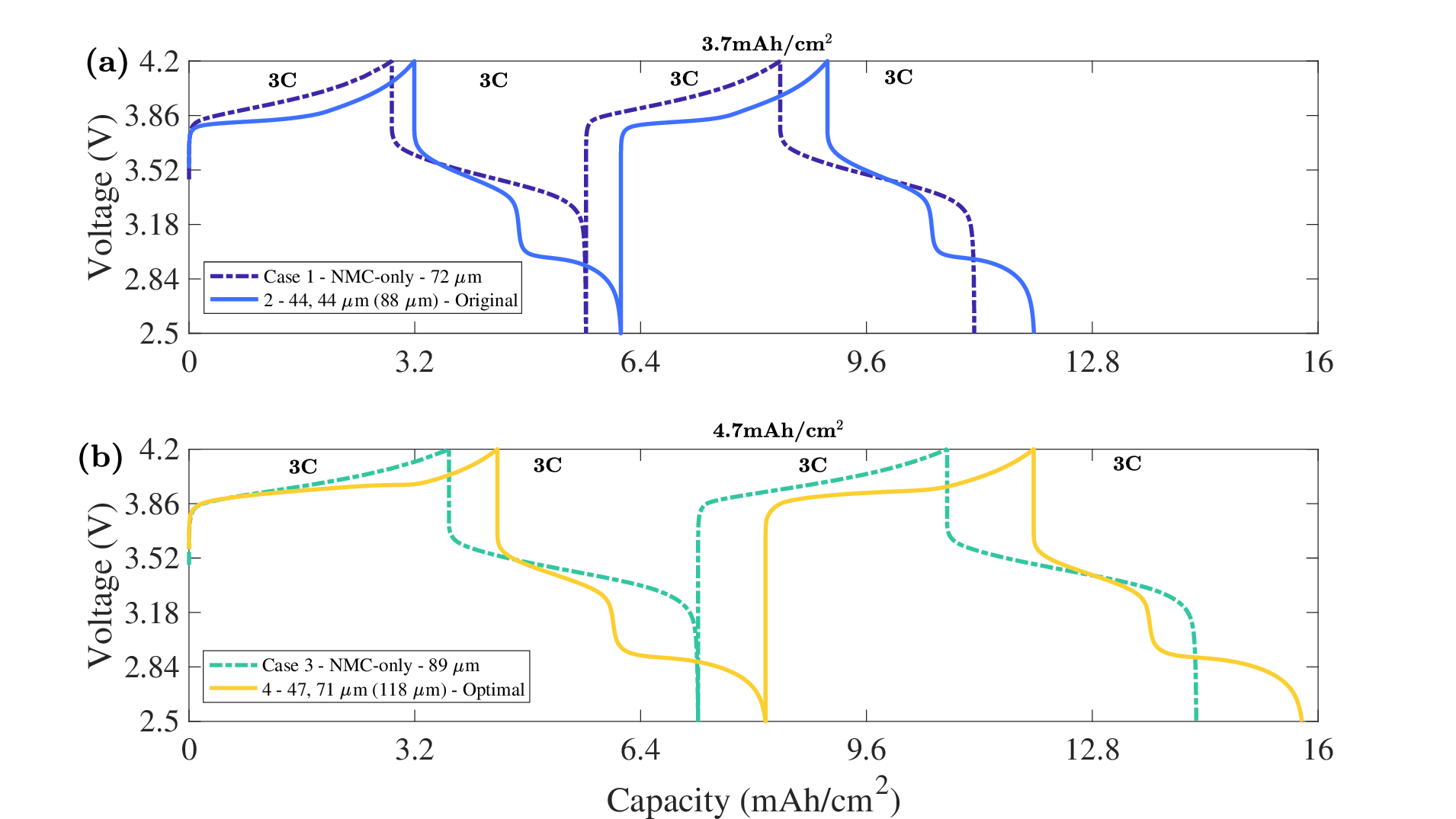}
  	\caption{ High C rate cycling at 3C including voltage as a function of capacity. The thinner original cells are shown in (a), while the optimal and equivalent specific capacity NMC-only electrode are in (b). For the bilayer in each cycle achieves capacities of 3.19, 2.93, 2.93, 2.93~mAh/cm$^2$ and 4.37, 3.80, 3.80, 3.80~mAh/cm$^2$ for the original and optimal, respectively.}
  	\label{fig:Cyclingfullsoc3c3c}
  \end{figure}

 We conduct short cycling tests to investigate the effect of high C rates on capacity. Two scenarios are investigated; \textit{(i)} a high C rate cycle at 3C to simulate usage, for example a power tool, along with \textit{(ii)} fast then slower cycling, at 3C charge, 0.1C discharge, then 3C charge, to simulate the usage, for example a laptop, mobile phone or EV. We utilise short cycles and will investigate the effects of longer cycles and degradation both within experiments and modelling in future work. The results are shown in Figs.~\ref{fig:Cyclingfullsoc3c3c} and \ref{fig:Cyclingfullsoc3c01c} including equivalent specific capacities at 0.05C grouped together in (a) for 3.7~mAh/cm$^2$ and (b) for 4.7~mAh/cm$^2$. The key point is that the capacity retention after the first charge is maintained in subsequent cycling and the positive effect of the increased capacity of the bilayer compounds with subsequent cycles. For 3C cycles in Fig.~\ref{fig:Cyclingfullsoc3c3c} the capacity of the subsequent discharge is lower, which is to be expected, as the retention is only 90\% after the first charge, and is starting at 90\% instead of 100\% charged. These trends are similar to most batteries at high C rates.

 In Fig.~\ref{fig:Cyclingfullsoc3c01c}, the 3C, 0.1C and 3C cycle is shown and the lower 0.1C in the subsequent discharge is sufficiently low enough to utilise all of the initial 90\% capacity retention, so is able to maintain the same capacity retention throughout. The largest capacity difference of the optimal bilayer compared to the NMC-only electrode out of Figs.~\ref{fig:Cyclingfullsoc3c3c} and \ref{fig:Cyclingfullsoc3c01c} is the 3C, 0.1C, 3C cycle in Fig.~\ref{fig:Cyclingfullsoc3c01c} (b), of 43.3\% collectively, or 14.4\% on average per cycle. To summarise, the new optimal bilayer is able to maintain high C rate charge on subsequent cycles, and has a significantly increased thickness and capacity. Compared to NMC-only electrodes of equivalent specific capacity at 0.05C, the gains are even higher with cycling for the optimal bilayer, which shows the benefits of creating an electrode using model-led design.

     \section*{Conclusion}
     
The multilayer Doyle-Fuller-Newman (M-DFN) model was used to optimise the design of NMC:LFP:CC lithium-ion battery electrodes for 3C charging. A cost function based on capacity retention was used when comparing various configurations to investigate capacity increase. The variable we were optimising was found to be the cathode active material mass (not the electrode length or sub-layer thickness ratio). The default and optimised case were validated with a 3D microscopic CT model, producing a close comparison. We show that $i)$ the sum of small gains in capacity retention produces larger gains of capacity and $ii)$ by making the electrolyte pathway slightly more porous and consequentially less tortuous, along with reducing the carbon binder amount, we were able to significantly increase the electrode thickness without sacrificing capacity retention. To achieve higher performance, we simultaneously update the electrochemical parameters accordingly:
     \begin{itemize}
     	\item Increase total electrode thickness,
     	\item Alter sub-layer thickness ratio,
     	\item Set electrolyte volume fraction the same at 0.3 for both sub-layers, 
     	\item Reduce LFP tortuosity,
     	\item Reduce carbon binder, which leads to increased surface area for reactions.
     \end{itemize}

The new optimal bilayer, which is 30~$\mu$m thicker than the default, charges between 0-90\% SOC in 18.6 minutes at a 3C rate, achieving 4.4~mAh/cm$^2$, and 15\% and 41\% higher than the NMC-only and LFP-only electrode, respectively. By creating a more homogeneous current response across the electrode through-thickness in terms of charge current density, the capacity was increased and electrolyte depletion reduced. Future work will focus on generalising the bilayer approach, and modelling and design optimisation frameworks, to full cells with bilayer cathodes and anodes.

\subsubsection*{Acknowledgements}
{\small 
The authors wish to acknowledge Patrick Grant, Yige Sun and Sam Wheeler for their contributions and discussions on the current and previous work\cite{Tredenick2024multilayer}. Funding for this work was provided by the Faraday Institution through project `Nextrode - next generation electrodes' (Grant numbers: FIRG015 and FIRG066). RD acknowledges funding by a UKIC Research Fellowship from the Royal Academy of Engineering.
}
\bibliographystyle{ieeetr}
\fontsize{8}{8}\rm
\bibliography{Battery_refs2}

 \clearpage

\subsection*{Appendix}

\renewcommand{\thefigure}{A\arabic{figure}}
\renewcommand{\theequation}{A\arabic{equation}}
\renewcommand{\thetable}{A\arabic{table}}
\setcounter{figure}{0} 
\setcounter{equation}{0} 
\setcounter{table}{0}

\begin{table} [h!]
	\footnotesize
	\centering
	\caption[Variables]{Variables of the M-DFN model, and $k$ $\in$ \{s,p\}, where $k$ is either the separator (s) or positive electrode (p).}	 	\label{VariablesModel}	
	\begin{tabular}{    p{1.2cm}  p{6.8cm}  p{1cm} } 
		\hline
		Parameter & Description  & Unit   \\ \hline 
		$a_{\scriptscriptstyle \text{p,i}}$	 &  Reaction surface area per volume  &  1/m   \\ 
		$c_{\scriptscriptstyle \text{e,k,i}}$ 	 &    Concentration of Li ions in the liquid electrolyte&  mol/m$^3$      \\ 
		$c_{\scriptscriptstyle \text{s,p,i}}$	 	 &  Concentration of Li ions in the solid particles   &   mol/m$^3$  \\ 
		$D_{\scriptscriptstyle \text{e,k,i}}$  	  	 &   Diffusivity function in the electrolyte  &   m$^2$/s    \\ 
		$D_{\scriptscriptstyle \text{e,bulk}}$  	  	 &   Bulk diffusivity function in the electrolyte &   m$^2$/s    \\ 
		$J_{\scriptscriptstyle \text{p,i}}$	 &  Butler-Volmer reaction current density  &  A/m$^2$   \\  
		$J_{\scriptscriptstyle \text{0,p,i}}$	 & Exchange current density &  A/m$^2$   \\ 
		$r_{\scriptscriptstyle \text{p,i}}$	 &   Particle radius direction  &   m   \\ 
		$t$ 	 &   Time  &  s    \\ 
		$U_{\scriptscriptstyle \text{p,i}}$ 	 &  Open circuit potential (OCP), equations (\ref{OCPNMC}) and (\ref{OCPLFP}).  &   V    \\ 
		$V$ 	 & Voltage    &   V   \\ 
		$x$	 &  Direction through cell   &  m    \\ 
		$\kappa_{\scriptscriptstyle \text{e,k,i}}$  	  	 &   Ionic conductivity function in the electrolyte  &   S/m  \\ 
		$\kappa_{\scriptscriptstyle \text{e,bulk}}$  	  	 & Bulk ionic conductivity function in the electrolyte &   S/m  \\ 
		$\eta_{\scriptscriptstyle \text{p,i}}$ &   Local overpotential  &   V   \\ 
		$\phi_{\scriptscriptstyle \text{s,p,i}}$  	 & Potential in the solid particles    &  V    \\ 
		$\phi_{\scriptscriptstyle \text{e,k,i}}$  	 & Potential in the liquid electrolyte    &  V    \\ 
		\hline 
	\end{tabular}
\end{table}

\begin{table} [h!]
	\footnotesize
	\centering
	\caption{Comparison of original and final optimal bilayer design parameters for the CT model including a cross sectional area of 9.604$\times 10 ^{-11}$~m$^2$ and reaction rate constant of $2\times 10^{-11}$ and $1\times 10^{-10}$~m$^{2.5}$s$^{-1}$mol$^{-0.5}$ for NMC and LFP, respectively. Additional parameters are as described in the previous work\cite{Tredenick2025CTvsDFN}.}	 	\label{optimaltable4}	  
	\begin{tabular}{lllll}
		\hline
		Parameter        &    &      Unit   &   Default  & Optimal   \\  	\hline
		NMC sublayer thickness     &  &  $\mu$m  & 44 &47    \\
		LFP sublayer thickness     &  &  $\mu$m  & 44 &71   \\
		Cathode electrode thickness     &  &  $\mu$m  & 88 &118  \\
		NMC thickness      &  &  \% &50.0 &39.8  \\ 
 	Applied current at 1C    &  &  mAh  &2.67$\times 10^{-6}$ & 3.21$\times 10^{-6}$  \\
		Specific areal capacity at 0.05C    &  &  mAh/cm$^2$  &3.712 &4.69    \\ \hline
		NMC  carbon binder domain volume fraction &   $\varepsilon_{\scriptscriptstyle \text{CBD,p1}}$     &   -&0.25& 0.25\\
		LFP  carbon binder domain volume fraction  &   $\varepsilon_{\scriptscriptstyle \text{CBD,p2}} $     & -&0.2  & 0.2\\
		Bruggeman tortuosity factor	&$b_{\scriptscriptstyle \text{p}}$&-& 1.5& 1.5  \\ 
		\hline
	\end{tabular}
\end{table}

 \section*{Open Circuit Potential}
 
 	\footnotesize
 The open circuit potential (OCP) functions\cite{Tredenick2024multilayer} of NMC622 and LFP using 0.02 C are shown below in equations (\ref{OCPNMC}) and (\ref{OCPLFP}). The NMC622 open circuit potential function was formulated by fitting to the experimental data, by first normalising the data from [0, 1] to [0.27, 0.92], then using a rational function:
 \begin{equation}
 	U_{\scriptscriptstyle \text{p1}} (x) =  \dfrac  { p1 \ x^4 + p2 \ x^3 + p3 \ x^2 + p4 \ x + p5 }{ x^5 + q1 \ x^4 + q2 \ x^3 + q3 \ x^2 + q4 \ x + q5},\label{OCPNMC}
 \end{equation}
 where $x$ is $c_{\scriptscriptstyle \text{p1}}/c^{\scriptscriptstyle \text{max}}_{\scriptscriptstyle \text{p1}}$, and
 \begin{multline*}
 	p1 =      -204.3   ,  \
 	p2 =      -166.6   ,  \
 	p3 =      -172.4  ,  \
 	p4 =       167.3   ,  \
 	p5 =       272.2  ,  \
 	q1 =      -158.1   ,  \
 	q2 =       221.4   ,  \
 	q3 =      -331.6  ,  \
 	q4 =       200.1   ,  \
 	q5 =       38.07   .
 \end{multline*}
 
 The open circuit potential of LFP is:    
 \begin{equation}
 	U_{\scriptscriptstyle \text{p2}} (y) =    3.413      +0.001 \left( \dfrac{1} { y} + \dfrac{1}{(y-1) }  \right) ,\label{OCPLFP} 
 \end{equation}
 where $y$ is $c_{\scriptscriptstyle \text{p2}}/c^{\scriptscriptstyle \text{max}}_{\scriptscriptstyle \text{p2}} $, and the value 3.413 is the mean of the plateau region of the OCP curve.

\begin{figure} [h!]
	\centering
	\includegraphics[width=0.3\textheight,keepaspectratio]{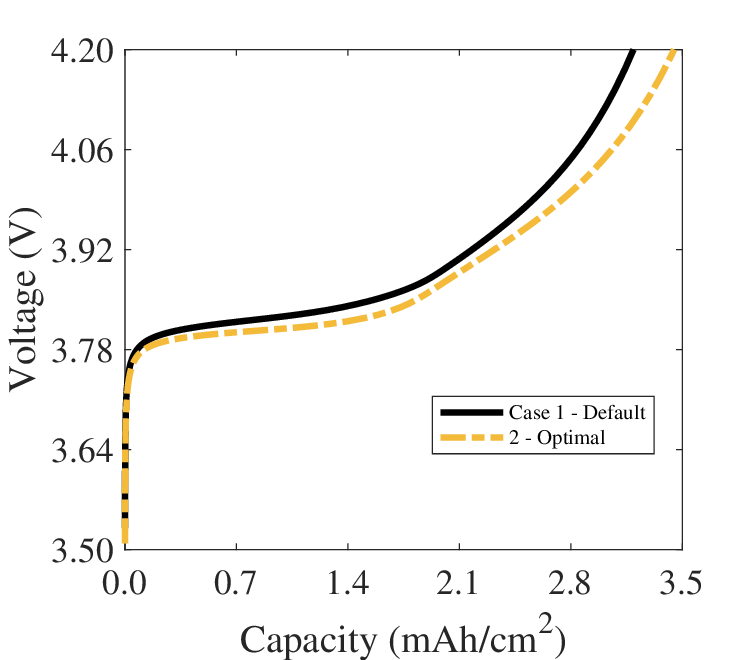}
	\caption{ The default bilayer compared to the new optimal during a 3C charge with parameters shown in Table \ref{optimaltable2}.}
	\label{fig:firstoptimal1}
\end{figure}

\begin{figure} [h!]
	\centering
	\includegraphics[width=0.3\textheight,keepaspectratio]{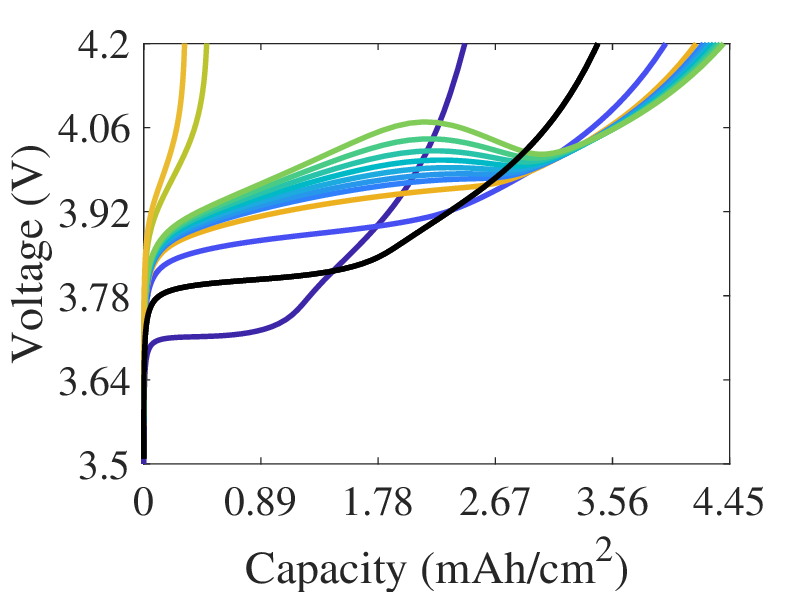} 
	\caption{ 3C charge with a range of cathode thicknesses. These are additional thicknesses, expanding upon Fig.~\ref{fig:compthtick2_3C}, shown in light blue and green.  }
	\label{fig:thickagain}
\end{figure}

\begin{figure} [h!]
	\centering
	\includegraphics[width=0.3\textheight,keepaspectratio]{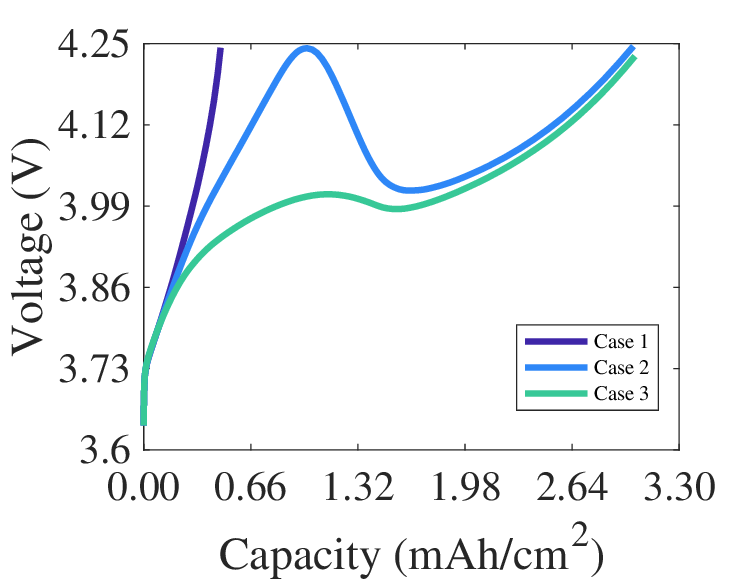} 
	\caption{ 3C charge with a range of LFP electrolyte volume fractions for the CT model. The cases have NMC electrolyte volume fractions of 0.2 and specific capacity at 0.05C of 3.7~mAh/cm$^2$. Cases 1, 2 and 3 have the LFP electrolyte volume fractions of 0.2, 0.4 and 0.6, respectively. Note here the cut-off voltage is set to 4.25~V instead of 4.2~V.
	}
	\label{fig:thickagainCT}
\end{figure}

 \begin{figure} [h!]
  	\centering
 	\includegraphics[width=0.4\textheight,keepaspectratio]{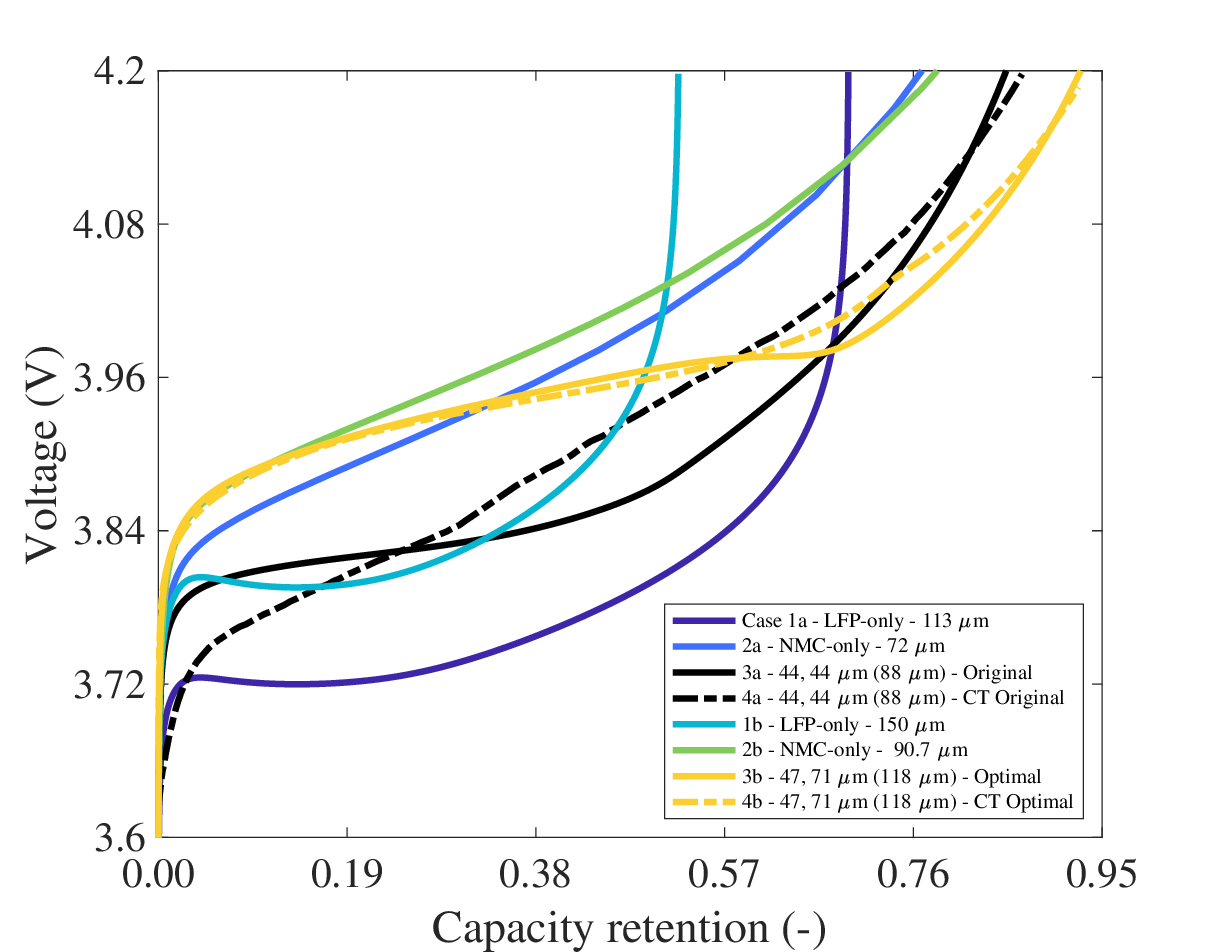}
 	\caption{Optimal of case 3b compared to the original bilayer of case 3a, during a 3C charge, expanding on Fig.~\ref{fig:optimal}. The NMC-only and LFP-only cases (1a, 2a, 1b, 2b) are included with equivalent specific capacity at 0.05C. Voltage is shown as a function of capacity retention. Cases 1a to 3a have capacity $\sim$3.7~mAh/cm$^2$ and cases 1b to 3b have capacity $\sim$4.7~mAh/cm$^2$. The thicker cases 1b to 3b include the candidate optimal design parameters in Table~\ref{optimaltable2} and \ref{optimaltable}.
 	}	
 	\label{fig:optimal3}
 \end{figure}

   \begin{figure} [h!]
	\centering
	\includegraphics[width=0.57\textheight,keepaspectratio]{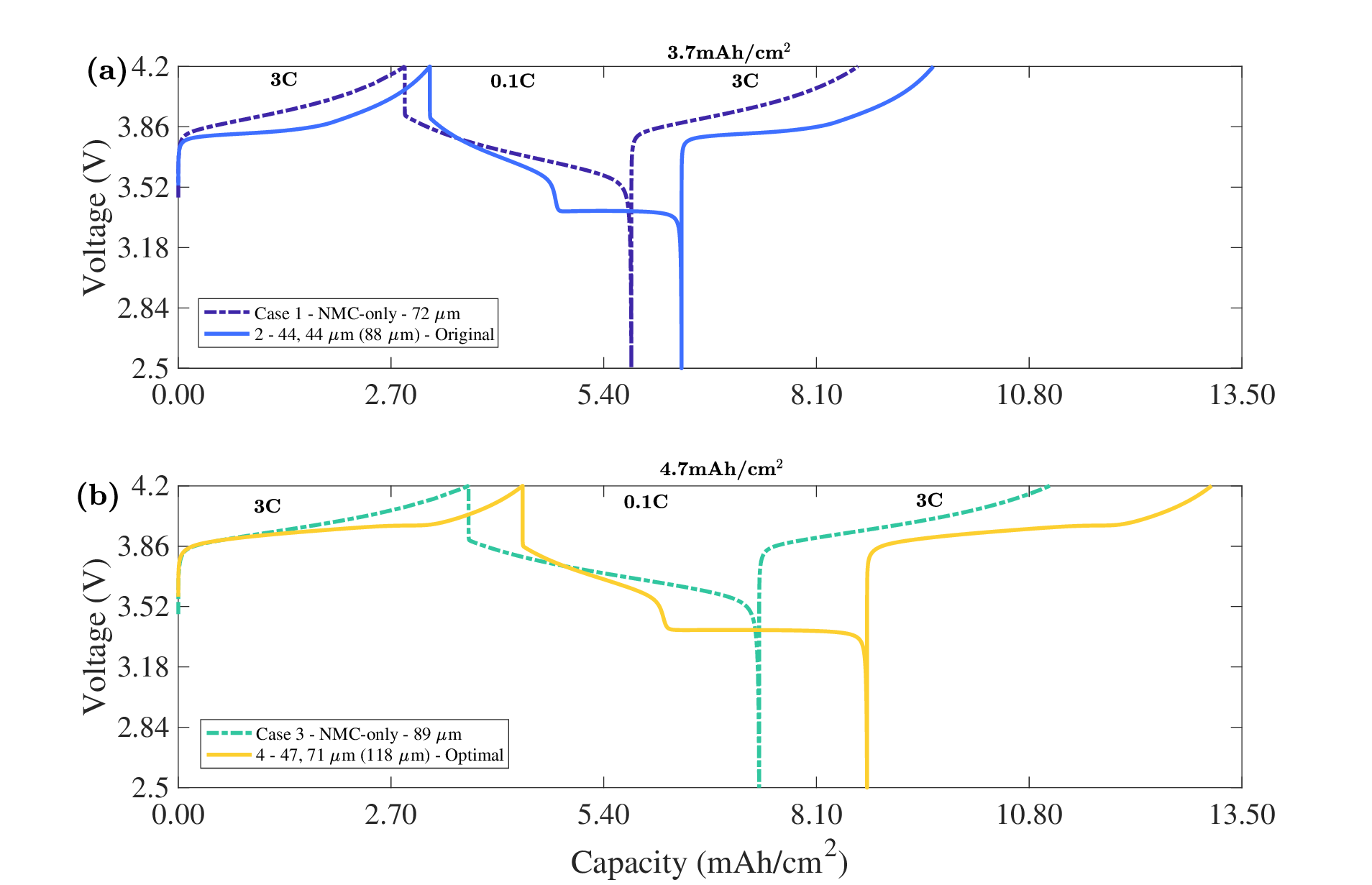}
	\caption{  Fast then slower C rate cycling of 3C charge, 0.1C discharge, 3C charge, including voltage as a function of capacity. The thinner original cells are shown in (a), while the optimal and equivalent specific capacity NMC-only electrode are in (b). For the bilayer in each cycle achieves capacities of 3.19, 3.19, 3.19~mAh/cm$^2$ and 4.37, 4.37, 4.37~mAh/cm$^2$ for the original and optimal, respectively.}
	\label{fig:Cyclingfullsoc3c01c}
\end{figure}

\end{document}